# Complexity problems associated with matrix rings, matrix semigroups and Rees matrix semigroups

Steve Seif    Željko Sokolović    Csaba Szabó [*]

July 29, 1997

## Abstract

Complexity problems associated with finite rings and finite semigroups, particularly semigroups of matrices over a field and the Rees matrix semigroups, are examined. Let $M_n\mathcal{F}$ be the ring of $n \times n$ matrices over the finite field $\mathcal{F}$ and let $T_n\mathcal{F}$ be the multiplicative semigroup of $n \times n$ matrices over the finite field $\mathcal{F}$. It is proved that for any finite field $\mathcal{F}$ and positive integer $n \geq 2$, the polynomial equivalence problem for the $T_n\,\mathcal{F}$ is co-NP-complete, thus POL-EQ$_\Sigma(M_n\,\mathcal{F})$ (POL-EQ$_\Sigma$ is a polynomial equivalence problem for which polynomials are presented as sums of monomials) is also co-NP-complete thereby resolving a problem of J. Lawrence and R. Willard and completing the description of POL-EQ$_\Sigma$ for the finite simple rings. In connection with our results on rings, we exhibit a large class of combinatorial Rees matrix semigroups whose polynomial equivalence problem is co-NP-complete. On the other hand, if $\mathbf{S}$ is a combinatorial Rees matrix semigroup with a totally balanced structure matrix $M$, then we prove that the polynomial equivalence problem for $\mathbf{S}$ is in P. Fully determining the complexity of the polynomial equivalence problem for combinatorial Rees matrix semigroups may be a difficult problem. We describe a connection between the polynomial equivalence problem for combinatorial Rees matrix semigroups and the retraction problem RET for bipartite graphs, a problem which computer scientists suspect may not admit a dichotomy into P and NP-complete problems (assuming P is not equal to NP).

For a combinatorial Rees matrix semigroup $\mathbf{S}$, a polynomial-time algorithm for the term-equivalence problem is given; the same is done for $\mathbf{S^1}$. In the last section, results for term-equivalence for Rees matrix semigroups are presented along with problems concerning term-equivalence for semigroups and general algebras. We exhibit an infinite sequence of finite subsemigroups $\mathbf{A_1} \leq \ldots \mathbf{A_k} \leq \mathbf{A_{k+1}} \ldots$ such that each semigroup in the sequence generates the same variety, for every positive integer $k$, POL-EQ$(\mathbf{A_{2k}})$ is in P, and POL-EQ$(\mathbf{A_{2k+1}})$ is co-NP-complete.

[*]All three authors were supported by the NSERC of Canada

# 1 Introduction

The complexity of equation solving and related problems in finite algebraic systems is a well-studied branch of complexity theory: Boolean algebras, vector spaces, rings, groups, and other algebraic systems give rise to complexity problems of interest (See [1], [2], [10],[11] for recent results of this type). We present results on complexity problems associated with semigroups which are then used to resolve complexity problems involving matrix rings.

The reader is assumed to be familiar with the basics of computational complexity. We use the standard terminology of computational complexity (see [7]), abbreviating polynomial-time complexity to P, non-deterministically polynomial time to NP, using $\equiv_P$ to denote the polynomial equivalence of two complexity problems and $\leq_P$ to denote the polynomial embedding of one complexity problem into another. Familiarity with standard complexity problems such as 3-SAT and GRAPH 3-COLORABILITY, as well as standard notions from graph theory, is also assumed.

A finite algebra $\mathbf{A} = (A; \mathcal{F})$ is a set, $A$, equipped with a finite indexed set $\mathcal{F} = \{f_1, \ldots, f_n\}$ of finitary *fundamental operations* on $\mathbf{A}$. In the sequel, all discussion is assumed to involve only finite algebras. If $f$ is a fundamental operation which involves $k$ variables, then we say that $f$ is *k-ary*. For more on general algebras see [3]. Complexity problems associated with $\mathbf{A}$ fall roughly into two categories. What is the computational complexity of determining whether an equation (or system of equations) has a solution? What is the computational complexity of determining whether two syntactic objects (eg. two polynomials) give rise to the same function of $\mathbf{A}$? We refer to complexity problems of these types as *algebra-complexity* problems. We will refer to *ring-complexity*, *semigroup-complexity* problems when our focus is restricted to rings or semigroups, respectively. In a sequel to this paper ([14]), we show that there is a correspondence between a semigroup-complexity problem, the polynomial equivalence problem (See POL-EQ, defined below) for the *combinatorial Rees matrix semigroups* and a currently researched graph-coloring complexity problem, the retraction problem RET, for bipartite graphs (see [6]).

## 1.1 Description of some algebra-complexity problems

We assume the reader is familiar with the definitions of semigroup, group, ring, field, and the basics of linear algebra. We review some algebra-complexity problems after presenting definitions of *terms* and *polynomials* for general algebras.

Let $(A; \mathcal{F})$ be a finite algebra. We define *terms* and *polynomials* over the algebra $\mathbf{A}$. Each variable $x_1, \ldots, x_n, \ldots$, is a term; if $f \in \mathcal{F}$ is a k-ary fundamental operation, and $g_1, \ldots, g_k$ are terms, then $f(g_1, \ldots, g_k)$ is a term. So terms are the syntactic objects formed from variables and fundamental operations. Polynomials are the syntactic objects formed from constants (elements of $A$), variables, and fundamental operations in the same way. If $p$ is a k-ary polynomial over $\mathbf{A}$, then $p$ induces a function from $A^k$ to $A$. Two polynomials $p$ and $q$ are said to be equal if they induce the same function on $A$.

We list some algebra-complexity complexity problems that have been taken up in the literature. TERM-EQ($\mathbf{A}$) has as an instance a pair of terms $\{p,q\}$ with the following question. Is $p$ equal to $q$? POL-EQ($\mathbf{A}$) has as instance a pair of polynomials $\{p,q\}$ with question: Is $p$ equal to $q$? POL-SAT($\mathbf{A}$) has as instance $\{p,b\}$, where $p$ is a polynomial and $b$ is an element of $A$ and the following question. With $p = p(x_1, \ldots, x_k)$, does there exist $(a_1, \ldots, a_k) \in A^k$ such that $p(a_1, \ldots, a_k) = b$ (i.e. does $p(x_1, \ldots, x_k) = b$ have a solution?). TERM-SAT has an instance $\{p,b\}$, where $p$ is a term and $b$ is an element of $A$ and asks if $p = b$ has a solution. The *length* of a term or a polynomial $p$ is defined to be the number of symbols in the expression (denoted by $length(p)$). We add the lengths of the various term and polynomial expression in an instance to arrive at the *size* of a given instance. Observe that TERM-EQ and POL-EQ are both in co-NP and that POL-SAT is in NP. Moreover, if POL-EQ($\mathbf{A}$) is in P, then TERM-EQ($\mathbf{A}$) is in P; if TERM-EQ($\mathbf{A}$) is co-NP-complete, then POL-EQ($\mathbf{A}$) is co-NP-complete. Later following [11] we describe ring-complexity problems for which instances involve terms or polynomials in a specified normal form.

## 1.2 Recent results on algebra-complexity problems for finite rings

In this subsection, we overview recent ring-complexity results of [1], [2], [10], and [11]. None of the results cited in this section play a role in the sequel. For further detail on finite ring theory see [9].

**Theorem 1.1 ([2], [10]).** *Let $\mathbf{R}$ be a finite ring. Then TERM-EQ($\mathbf{R}$) is in P if and only if $\mathbf{R}$ is nilpotent; otherwise, TERM-EQ($\mathbf{R}$) is co-NP-complete.*

An examination of the proof of [10] reveals that POL-EQ($\mathbf{R}$) is in P if $\mathbf{R}$ is nilpotent. Also TERM-EQ($\mathbf{R}$) is co-NP-complete for non-nilpotent finite rings (which is proved in [2]) implies POL-EQ($\mathbf{R}$) is co-NP-complete for these rings. Thus POL-EQ for finite rings is also resolved by [2], [10]: for a finite ring $\mathbf{R}$, POL-EQ($\mathbf{R}$) is in P if and only if $\mathbf{R}$ is nilpotent; otherwise, POL-EQ($\mathbf{R}$) is co-NP-complete.

A *monomial* over the ring $\mathbf{R}$ is a product of non-commuting variables and elements of $\mathbf{R}$ (so "constants" are allowable symbols in a monomial). In [11] the authors consider algebra-complexity problems for rings, requiring that in a given instance, occurences of polynomials and terms be sums of monomials. In [11], for each of the ring-complexity problems defined above a related problem is defined, one in which terms and polynomials in a given instance are required to be sums of monomials. We add the subscript "$\Sigma$" to indicate instances involve sums of monomials (eg. TERM-EQ$_\Sigma(\mathbf{R})$, POL-EQ$_\Sigma(\mathbf{R})$, and POL-SAT$_\Sigma(\mathbf{R})$).

Restricting instances to the terms and polynomials that are sums of monomials has an effect on ring-complexity questions. For example if we consider the two element ring $\mathbf{Z_2}$, it is straight-forward consequence of the fact that 3-SAT is NP-complete (see [7]) that TERM-EQ($\mathbf{Z_2}$) is co-NP-complete; whereas, TERM-EQ$_\Sigma(\mathbf{Z_2})$ is in P, as can be easily verified.

**Theorem 1.2.** *[11] Let $\mathbf{R}$ be a finite ring with Jacobson radical $\mathcal{J}$.*

1. *If $\mathbf{R}/\mathcal{J}(\mathbf{R})$ is commutative, then TER-EQ$_\Sigma(\mathbf{R})$ is in P.*

2. *If $\mathbf{R}$ is a finite simple matrix ring whose invertible elements form a non-solvable group, then TER-EQ$_\Sigma(\mathbf{R})$ is co-NP-complete.*

## 1.3 Statement of the main results

We begin with some notation.

- Let $M_n\ \mathcal{F}$ be the ring of $n$ by $n$ matrices over the field $\mathcal{F}$.

- Let $T_n\ \mathcal{F}$ denote the multiplicative semigroup of $n \times n$ matrices over $\mathcal{F}$.

- Let $L_n\ \mathcal{F}$ be the subsemigroup of $T_n\ \mathcal{F}$ consisting of matrices of rank 1, along with the 0-matrix.

**Theorem 1.3.**  1. *The computational complexity of POL-EQ($L_n\ \mathcal{F}$) is co-NP-complete.*

2. *The computational complexity of POL-SAT($L_n\ \mathcal{F}$) is NP-complete.*

**Theorem 1.4.**  1. *The computational complexity of POL-EQ($T_n\ \mathcal{F}$) is co-NP-complete.*

2. *The computational complexity of POL-SAT($T_n\ \mathcal{F}$) is NP-complete.*

If a problem concerning monomials over the ring is intractable, then the corresponding problem for the ring is also intractable; hence, the next theorem is a consequence of Theorem 1.4.

**Theorem 1.5.**  1. *The computational complexity of POL-EQ$_\Sigma(M_n\ \mathcal{F})$ is co-NP-complete.*

2. *The computational complexity of POL-SAT$_\Sigma(M_n\ \mathcal{F})$ is NP-complete.*

Let $\mathbf{G}$ be a finite group. In [12], it is proved that if $\mathbf{G}$ is not solvable, then TER-EQ($\mathbf{G}$) is co-NP-complete. In [11] the authors polynomially embed the term equivalence problem for $GL_n\ \mathcal{F}$ into the term equivalence problem for $M_n\ \mathcal{F}$ and in so doing prove that the TER-EQ$_\Sigma M_n\ (\mathcal{F})$ is co-NP-complete if $GL_n\ \mathcal{F}$ is non-solvable.

Recall that a finite ring is simple if and only if it is isomorphic to the full matrix ring of $n \times n$ matrices over a field. For a field $\mathcal{F}$, it is straightforward to prove that POL-SAT$_\Sigma\ (\mathcal{F})$ is in P (or use Theorem 1.2, item 2). Thus Theorem 1.5 completes the description of POL-EQ$_\Sigma$ and POL-SAT$_\Sigma$ for the simple rings. Recall that a ring $\mathbf{R}$ is said to be *semi-simple* if its Jacobson radical $\mathcal{J}(\mathcal{R})$ is trivial and that a finite ring is semi-simple if and only if it is the direct product of simple rings.

**Corollary 1.6.** *Let $\mathbf{R}$ be a finite semi-simple ring.*

1. *POL-EQ$_\Sigma$(**R**) is in P if and only if **R** is commutative. Otherwise, POL-EQ$_\Sigma$(**R**) is co-NP-complete.*

2. *POL-SAT$_\Sigma$(**R**) is in P if and only if **R** is commutative. Otherwise, POL-SAT$_\Sigma$(**R**) is NP-complete.*

### 1.3.1 Semigroup-complexity for combinatorial Rees matrix semigroups

Later in the paper we provide implementable polynomial time algorithms for the term-equivalence problem for combinatorial Rees matrix semigroups (See Lemma 4.15). The reader unfamiliar with Rees matrix semigroups may want to consult the first two pages of the next subsection.

**Theorem 1.7.** *Let **S** be a combinatorial finite Rees matrix semigroup. Then TERM-EQ(**S**) is in P.*

Toward the long-range goal of classfying TERM-EQ problems for semigroups, for a Rees matrix semigroup **S**, one would like to determine conditions under which an extension **T** satisfies TERM-EQ(**T**)$\equiv_P$ TERM-EQ(**S**).

**Theorem 1.8.** *Let **S** be a combinatorial Rees matrix semigroup. Then TERM-EQ(**S**$^1$) is in P.*

Let **S** be a Rees matrix semigroup. For the next problem we assume a background in semigroup theory. Let $\Omega(\mathbf{S})$ be the translational hull of **S**. Of course, **S**$^1$ is a subsemigroup of $\Omega(\mathbf{S})$. Theorem 1.8 motivates the following problem.

**Problem 1.9.** *Describe the (combinatorial) Rees matrix semigroups **S** such that TERM-EQ($\Omega(\mathbf{S})$) is in P.*

Resolving the term equivalence problem for translational hulls of Rees matrix semigroups may contribute greatly toward the larger goal of providing a complete description of the term-equivalence problem for semigroups. In Lemma 5.5, we compute TERM-EQ($\Omega(\mathbf{S})$) for certain Rees matrix semigroups **S**.

It is well-known that the semigroups $\mathcal{L}_n Z_2$ are isomorphic to combinatorial Rees matrix semigroups (see Lemma 2.2 and its proof), so Theorem 1.3 provides examples of combinatorial Rees matrix semigroups with co-NP-complete polynomial-equivalence problems. The next theorem shows that there exist large classes of combinatorial Rees matrix semigroups with polynomial-equivalence problems that are in P.

A 0-1 matrix $M$ is said to be *totally balanced* if $M$ satisfies for all $\alpha, \beta \in \Lambda$, and all $i, j \in I$ with $i \neq j$, $\alpha \neq \beta$, $M(\alpha, i) = 1 = M(\alpha, j) = M(\beta, i) = 1$, then $M(\beta, j) = 1$.

**Theorem 1.10.** *Let $M$ be a totally balanced matrix. Then POL-EQ(**S**$_\mathbf{M}$) is in P.*

This completes the description of our main results. In Section 2 we review basic definitions and results concerning Rees matrix semigroups and begin the discussion of semigroup-complexity problems for Rees matrix semigroups. The proofs for the results described above are contained in Section 3 and Section 4. In Section 5 we prove some term equivalence results for general Rees matrix semigroups and pose additional problems.

## 2  Background information on Rees matrix semigroups

We begin with notation and definitions from semigroup theory. See [5] for further background on semigroup theory. Let **S** be a semigroup. If **S** is a subsemigroup of a semigroup **T**, we write $\mathbf{S} \leq \mathbf{T}$ and say that **T** is an *extension* of **S**; **S** is an ideal of **T** if for all $t \in T$ and all $s \in S$, $st, ts \in S$ and we say that **T** is an *ideal extension* of **S**. If $\mathbf{S} \leq \mathbf{T}$ and **S** satisfies the group axioms, then we say that **S** is a *subgroup* of **T**. An element 0 satisfying for all $s \in S$, $s0 = 0 = 0s$ is said to be a *zero* of **S**; an element 1 of $S$ satisfying for all $s \in S$ $1s = s = s1$ is said to be a *identity* of **S** (0 will always be used to denote the zero of a semigroup; 1 will be used to denote the identity). By $\mathbf{S^0}$ we mean the extension of **S** with underlying set $S \cup \{0\}$ where 0 is a zero. In like manner $\mathbf{S^1}$ will denote the extension of **S** with underlying set $S \cup \{1\}$ where 1 is an identity. For a semigroup **S** with a 0, an element $s \in S$ will be said to be a *0-divisor* if there exists $t \in S \setminus \{0\}$ such that either $st = 0$ or $ts = 0$.

The ideals of **S** form a partially ordered set under set inclusion. For subsets $U, V$ of the semigroup **S**, we let $UV = \{uv : u \in U, v \in V\}$. Note that every finite semigroup has a unique minimal ideal since the product of any finite set of ideals of **S** is non-empty and is an ideal. The next piece of notation is non-standard. If **S** has a 0 we say that an ideal $I$ of **S** is *0-minimal* if $I$ is minimal among all ideals properly containing the one element ideal $\{0\}$; $I$ is said to be *completely 0-simple* if $I$ is 0-minimal and satisfies Lemma 2.1, item 1.

Let **S** be a semigroup. Recall that a binary relation $\theta$ on $S$ is *compatible* if for all $(a,b) \in \theta$ and all $s, t \in S^1$, $(sat, sbt) \in \theta$. A compatible equivalence relation $\theta$ on $S$ is said to be a *congruence* of **S**. For a congruence $\theta$ of **S**, we let $\mathbf{S}/\theta$ denotes the quotient semigroup with underlying set $S/\theta$, the equivalence classes of $\theta$. For $s \in S$ we let $\overline{s}$ denote $s/\theta$. If $\theta$ is a congruence and $p$ is a polynomial, then $p$ induces a well-defined function on $\mathbf{S}/\theta$ which we denote by $\overline{p}$. In fact $\overline{p}$ is a polynomial over $\mathbf{S}/\theta$.

We give the preliminaries for the definition of a Rees matrix semigroup. Let **G** be a group. An $m \times n$ matrix $M$ with entries in $G \cup \{0\}$ is said to be a matrix over **G**. If every row and every column of $M$ has at least one element which is not 0, then $M$ is said to be a *regular* matrix. We let $\Lambda$ be an index set for the rows of $M$; we let $I$ be an index set for the columns of $M$. If $\lambda \in \Lambda$ and $i \in I$, then $M(\lambda, i)$ denotes the $(\lambda, i)$ entry of $M$. For any group **G**, we let 1 denote the identity of **G**. A *0-1 matrix* denotes a matrix over the trivial group.

For the remainder of the section, we let $M$ be an $m \times n$ regular matrix over a group **G**. We define the Rees matrix semigroup $\mathbf{S} = \mathcal{M}(G, M)$ as follows. The underlying set consists of 0 and all triples of the form $[i, g, \lambda]$ where $i \in I$, $g \in G$ and $\lambda \in \Lambda$, along with 0. An associative multiplication is given by the following rule: If $M(\lambda, j) \neq 0$, then $[i, g, \lambda][j, h, \gamma] = [i, gM(\lambda, j)h, \gamma]$. If $M(\lambda, j) = 0$, then $[i, g, \lambda][j, h, \gamma] = 0$.

Our notation is slightly non-standard. More standard is $\mathcal{M}^0(G, \Lambda, I, M)$ for the Rees matrix semigroup defined above. In the sequel, all matrices $M$ over a group **G** will be assumed to be regular.

For $\mathcal{M}(G, M)$, we say that **G** is the *structure* group of $\mathcal{M}(G, M)$; $M$ is the *structure matrix* of $\mathcal{M}(G, M)$. Though it plays no role in the sequel, it is not difficult to verify that $\mathcal{M}(G, M)$ has up to isomorphism a unique maximal subgroup which is isomorphic to its stucture group **G**. A Rees matrix semigroup $\mathcal{M}(G; M)$ is said to be *combinatorial* (or *group-free*) if its stucture group **G** is trivial. All results here directly involve or depend upon resolving semigroup-complexity problems for combinatorial Rees matrix semigroups with 0-divisors.

Let $\mathbf{S} = \mathcal{M}(G; M)$ be a Rees matrix semigroup. We define an equivalence relation $\mathcal{H}$ on $S$. For $[i, g, \lambda], [j, h, \beta] \in S$, $([i, g, \lambda], [j, h, \beta]) \in \mathcal{H}$ if $i = j$ and $\beta = \lambda$; the $\mathcal{H}$ class of 0 is $\{0\}$. We define an $m \times n$ 0-1 matrix $\overline{M}$ as follows. For $\overline{M}(\lambda, i) = 1$, if $M(\lambda, i) \in G$; $\overline{M}(\lambda, i) = 0$, if $M(\lambda, i) = 0$. It is straightforward to verify that $\mathcal{H}$ is a congruence of **S** and that $\mathbf{S}/\mathcal{H}$ is isomorphic to $\mathcal{M}(\{1\}, \overline{M})$. A combinatorial Rees matrix semigroup is entirely determined by its 0-1 Rees matrix $M$. We shorten $\mathcal{M}(\{1\}, M)$ to $\mathbf{S_M}$. In $\mathbf{S_M}$, we abbreviate $[i, 1, \lambda]$ to $[i, \lambda]$ ($i \in I$, $\lambda \in \Lambda$).

Lemma 2.1 consists of several well-known facts that are used repeatedly in the sequel. We suggest that the reader verify items 1-4. All have short proofs. For items 5-6, it is not difficult to guess the explicit isomorphisms.

For $\alpha \in \Lambda$, we let $M_\alpha$ denote the $\alpha$ row of $M$. For $i \in I$, let $M^{(i)}$ denote the i-th column of $M$. Two rows of $M$ say $M_\alpha$ and $M_\beta$, $\alpha, \beta \in \Lambda$, are said to *proportional* if there exists $g \in G$ such that for all $i \in I$, $gM(\alpha, i) = M(\beta, i)$. Two columns of $M$ say $M^{(i)}$ and $M^{(j)}$, $i, j \in I$, are said to *proportional* if there exists $g \in G$ such that for all $\lambda \in \Lambda$, $M(\lambda, i)g = M(\lambda, j)$.

Let $\mathbf{S} = \mathcal{M}(G; M)$ be a Rees matrix semigroup. For $s = [i, g, \lambda] \in S$ ($i \in I$, $g \in G$, and $\lambda \in \Lambda$), let $s_1 = i$ and let $s_2 = \lambda$ (we are ignoring the middle component for reasons that will become clear).

**Lemma 2.1.** *Let $\mathbf{S} = \mathcal{M}(G; M)$ be a Rees matrix semigroup with elements $a(1), \ldots, a(k) \in S$.*

1. *For $s, t \in S$, $s \neq 0$ and $t \neq 0$, there exist $u, v \in S$ such that $usv = t$.*

2. *$a(1) \ldots a(k) = 0$ if and only if there exists $j$, $1 \leq j \leq k$, such that $a(j-1)a(j) = 0$.*

3. *$a(1) \ldots a(k) = 0$ if and only if $\overline{a(1)} \ldots \overline{a(k)} = 0$ in $S/\mathcal{H}$. In particular, if $p$ is a polynomial, $p$ is identically 0 if and only if $\overline{p}$ is identically 0.*

4. *If $a(1) \ldots a(k) \neq 0$, then $(a(1) \ldots a(k))_1 = a(1)_1$ and $(a(1) \ldots a(k))_2 = a(k)_2$. In particular, if **S** is combinatorial, then $a(1) \ldots a(k) \neq 0$ implies $a(1) \ldots a(k) = [a(1)_1, a(k)_2]$.*

5. *Let $M'$ result from a sequence of permutations of the rows and columns of $M$. Then $\mathcal{M}(G; M')$ is isomorphic to $\mathcal{M}(G; M)$.*

6. *For $\lambda \in \Lambda$, replace $M_\lambda$ by a row which is proportional to $M_\lambda$ and denote the resulting $m \times n$ matrix by $M'$. Then $\mathcal{M}(G; M')$ is isomorphic to $\mathcal{M}(G; M)$. In like manner, if $i \in I$ if $M'$ results by replacing $M^{(i)}$ by a column which is proportional to $M^{(i)}$, then $\mathcal{M}(G; M')$ is isomorphic to $\mathcal{M}(G; M)$.*

The remarks contained in this paragraph play no role in the sequel. The property of Rees matrix semigroups described by Lemma 2.1, item 1 is of importance in semigroup theory. A semigroup **S** with 0 is said to be *completely 0-simple* if it satisfies Lemma 2.1, item 1; thus, **S** is completely 0-simple if **S** is a completely 0-minimal ideal of itself. The Rees Matrix Theorem states that every completely 0-simple semigroup is isomorphic to a Rees matrix semigroup.

## 2.1 Matrix rings over fields and Rees matrix semigroups

The results described in this section are well-known. For a matrix $M$, let $M^T$ denote the transpose of $M$. For a n-vector $v$ ($n$ is a positive integer) over the field $\mathcal{F}$, as is convenient in the sketch below, we regard $v$ as a $n \times 1$ matrix. For a field $\mathcal{F}$, let $\mathcal{F}^*$ denote the non-0 zero elements of $\mathcal{F}$, the multiplicative group of units of $\mathcal{F}$.

**Lemma 2.2.** *The semigroup $L_n \mathcal{F}$ is isomorphic to a Rees matrix semigroup.*

*Proof.* We provide a sketch of the proof, leaving some details to the reader. Let **M** denote $M_n \mathcal{F}$, let $\mathbf{L} = L_n \mathcal{F}$, $V = \mathcal{F}^n$ and let $I$ index the one-dimensional subspaces of $V = \mathcal{F}^n$. Let $|I| = r$ ($r = (|\mathcal{F}|^n-1)/(|\mathcal{F}|-1)$). From each 1-dimensional subspace $U_i$ ($i \in I$), choose a non-0 vector $v_i \in U_i$. Consider the ordered set $\{v_1, \ldots, v_r\}$. Let $A \in L$ be a rank-1 matrix. Observe that the kernel of $A$ is a subspace of codimension 1; hence, it is the orthogonal complement of a unique 1-dimensional subspace of $V$. Thus the kernels of the rank-1 matrices in $L$ can also be indexed by $I$. Let $A$ be a rank-1 matrix with range spanned by $v_i$ ($i \in I$) and whose kernel is the orthogonal complement of $v_\lambda$ ($\lambda \in I$). Verify that there exist a unique $k \in \mathcal{F}^*$ such that $A = (kv_i)v_\lambda^T = v_i(kv_\lambda^T)$. Map $A$ to the triple $[i, k, \lambda]$. Map the 0-matrix of $L$ to 0. Define an $r \times r$ matrix $T$ as follows. For $i, \lambda \in I$, let $T(\lambda, i) = v_\lambda^T v_i \in \mathcal{F}$.

The hardest part of this sketch will be to verify that the map described above is a semigroup isomorphism from $L$ to the Rees matrix semigroup $\mathcal{M}(\mathcal{F}^*; T)$. This completes our sketch. □

Let $m, n$ be positive integers. We let $J_{m \times n}$ be the $m \times n$ "all 1's" matrix: for $1 \leq i \leq m$, $1 \leq j \leq n$, $J_{m \times n}(i, j) = 1$. Let $I_k$ denote the $k \times k$ identity matrix and let $H_k = J_{k \times k} - I_k$.

**Lemma 2.3.** *For any field $\mathcal{F}$, the semigroup $(L_2\ \mathcal{F})/\mathcal{H}$ is isomorphic to $\mathbf{S_{H_n}}$, where $n = |\mathcal{F}| + 1$ is the number of 1-dimensional subspaces of the vector space $\mathcal{F}^2$. In particular, the semigroup $L_2\ \mathbf{Z_2}$ is isomorphic to $\mathbf{S_{H_3}}$.*

*Proof.* By Lemma 2.2, we have that $L_2\ \mathcal{F}$ is isomorphic to a Rees matrix semigroup $\mathcal{M}\ (\mathcal{F}^*, T)$.

We examine the construction of the Rees matrix $T$ in Lemma 2.2. Suppose $\{v_1, \ldots, v_r\}$ is the set defined in Lemma 2.2, an indexed set of non-zero vectors, one from each 1-dimensional subspace of $V$. For $\lambda \in I$, observe since the dimension of the vector space is 2, $v_\lambda$ is orthogonal to exactly one vector in $\{v_1, \ldots, v_r\}$. Hence each row of the square matrix $T$ has exactly one 0. By basic linear algebra, if $\alpha \neq \beta \in I$, then for all $i \in I$, if $v_\alpha^T v_i = 0$, then $v_\beta^T v_i \neq 0$. Thus no pair of distinct rows of $T$ has a 0 in the same column. It follows that $(L_2\ \mathcal{F})/\mathcal{H}$ is isomorphic to $\mathbf{S_{H_r}}$, where $r$ is the number of 1-dimensional subspaces of $\mathcal{F}^2$ (Arriving at the matrix $H_r$ may require column and rows exchanges. By Lemma 2.1, item 5, column and row exchanges leave the isomorphism class of the associated Rees matrix semigroup unchanged).

Since the number of 1-dimensional subspaces of $Z_2^2$ is 3, the last sentence of the lemma follows from the fact that $L_2\ Z_2$ is already combinatorial. $\square$

## 2.2 Some combinatorics for term and polynomial equivalence problems of Rees matrix semigroups

Our main technique for resolving semigroup-complexity problems involves analyzing the 0's of a polynomial $p$ by associating certain graphs to $p$. Let $\mathbf{S}$ be a semigroup with 0 and let $p$ be a k-ary polynomial over a semigroup $\mathbf{S}$. The set $\{(a_1, \ldots, a_k) : p(a_1, \ldots, a_k) = 0\}$ will be denoted $Z(p)$ and called a *Z-set*.

For polynomials $p, q$ over $\mathbf{S}$ we say that $p$ *syntactically divides* $q$ if there exists $s, t$ polynomials over $\mathbf{S}$ such that $q = spt$. For a polynomial $p$ over $\mathbf{S}$ we define the directed graph $G(p) = (V, E)$. Let $U$ denote the variables that occur in $p$ and let $C$ denote the constants that occur in $p$. We let the vertex set $V = U \cup C$ and for $c, d \in V$, $c, d$ is an edge (denoted $c, d \in E$) if $cd$ syntactically divides $p$.

The next lemma gives an indication that problems involving Z-sets can be solved in a combinatorial framework. In this section and later in the paper, we will show that the complexity of equivalence problems for certain semigroups is significantly determined by complexity problems involving Z-sets.

**Lemma 2.4.** *Let $\mathbf{S}$ be a Rees matrix semigroup and let $p, q$ be two polynomials over $\mathbf{S}$. If $G(p) = G(q)$, then $Z(p) = Z(q)$.*

*Proof.* By Lemma 2.1, item 2, a product in $\mathbf{S}$ is 0 if and only if the product of two adjacent factors is 0. Assume $G(p) = G(q)$. Since $G(p)\ (= G(q))$ encodes adjacencies of constants and variables occuring in $p$ and $q$, it follows that for any evaluation $e$ of the variables occuring in $p$ and $q$, $e(p) = 0$ if and only if $e(q) = 0$. $\square$

In general, the converse does not hold. For example, consider $I_2$, the $2 \times 2$ identity matrix and the Rees matrix semigroup $\mathbf{S_{I_2}}$. Verify that $p = x^2y^2 = y^2x^2 = q$ is an identity of $\mathbf{S_{I_2}}$. In particular, $Z(p) = Z(q)$. On the other hand, $G(p) \neq G(q)$. As we will show later, the converse does hold for the terms of an important class of Rees matrix semigroups.

We introduce a semigroup complexity problems for a semigroup $\mathbf{S}$ with a 0. Let POL=0($\mathbf{S}$) have as instance a polynomial $p$. We ask if $p$ is identically 0.

Let $\mathbf{S}$ be a semigroup, let $p$ be a polynomial over $\mathbf{S}$ with variables $U$ occuring in $p$. A function $e : U \to S$ is said to be an *evaluation*. A partial function $e : U \to S$ is said to be a *partial evaluation*. If $e$ is an evaluation, by $e(p)$ we mean the polynomial over $\mathbf{S}$ that results by substituting $e(x_i)$ for $x_i$ for each $x_i$ in the domain of $e$. For convenience, we include $S$ in the domain of any partial evaluation by defining $e(s) = s$, for all $s \in S$. For a polynomial $p$ over a semigroup, let $p_l$ denote the left-most symbol of $p$ and let $p_r$ denote the right-most symbol of $p$.

**Lemma 2.5.** *Assume $M$ be a $m \times n$ regular 0-1 matrix, let $\mathbf{G}$ be a group, and let $\mathbf{S} = \mathcal{M}(G; M)$. Assume that $\mathbf{T}$ is an ideal extension of $\mathbf{S}$ and $\mathbf{S}$ is a completely 0-minimal ideal of $\mathbf{T}$. Then the following hold.*

1. *POL=0($\overline{\mathbf{S}}$)$\equiv_P$POL=0($\mathbf{S}$).*

2. *POL=0($\overline{\mathbf{S}}$) $\leq_P$ POL-EQ($\mathbf{T}$).*

3. *If POL=0($\overline{\mathbf{S}}$) is co-NP-complete, then POL-SAT($\mathbf{T}$) is NP-complete.*

4. *POL=0($\overline{\mathbf{S}}$) is in P if and only if POL-SAT($\overline{\mathbf{S}}$) is in P.*

*Proof.* Item 1 follows immediately form Lemma 2.1, item 3, and the fact that the mapping $\{p\} \to \{\overline{p}\}$ is linear in size. By item 1, we can use $\mathbf{S}$ in place of $\overline{\mathbf{S}}$ in the proof of items 2 and 3.

Let $\{p\}$ be an instance of POL=0($\mathbf{S}$). Let $U = \{x_1, \ldots, x_n\}$ be the variables involved in $p$. Let $s \in S$ be a non-0 element. For $x_i \in U$, let $\rho(x_i) = x_i s y_i$. For each $x_i \in U$, replace each occurence of $x_i$ in $p$ by $\rho(x_i)$. Denote the polynomial that results by $\rho(p)$. We use $\rho(p)$ in the proofs of items 2,3 and 4.

The fact that $\mathbf{S}$ is a completely 0-minimal ideal of $\mathbf{T}$ is exactly what is needed to show that $p$ is identically 0 over $\mathbf{S}$ if and only if $\rho(p)$ is identically 0 over $\mathbf{T}$. To complete the proof of item 2, use that the mapping of instances $\{p\} \to \{\rho(p), 0\}$ from POL=0($\mathbf{S}$) to POL-EQ($\mathbf{T}$) is linear with respect to the size of instance.

We prove item 3. Assume $y_1$ and $y_2$ are variables not involved in $p$, a polynomial over $\mathbf{S}$. Let $p' = y_1 p y_2$. The regularity of $M$ implies that $p'$ is identically 0 over $\mathbf{S}$ if and only if $p$ is identically 0 over $\mathbf{S}$. Moreover for any $b = [i, g, \lambda] \in S \setminus \{0\}$, we claim that $p' = b$ has a solution in $\mathbf{S}$ if and only if $p$ is **not** identically 0 over $\mathbf{S}$. If $p$ is not identically 0, there exists an evaluation of $e$ of the variables of $p$ such that $e(p) \neq 0$. Let $e(p) = [r, h, \alpha]$ ($r \in I$, $h \in G$, $\alpha \in \Lambda$). Since $M$ is regular, there exists $\beta \in \Lambda$ and $s \in I$ such that $M(\beta, r) \neq 0$ and $M(\alpha, s) \neq 0$. Using the fact that $\mathbf{G}$ is a group, it follows that there exist $j, k \in G$ such that $[i, j, \beta][r, h, \alpha][s, k, \lambda] = [i, g, \lambda] = b$. This completes one direction of the claim; the other direction is obvious. Using the fact that **S** is a completely 0-minimal ideal of **T**, it follows that for $b \in S \setminus \{0\}$, $\rho(p') = b$ has a solution in **T** if and only if $p' = b$ has a solution in **S**. By putting this sentence together with the claim preceding it, we have $p$ is identically 0 if and only for any $b \in S \setminus \{0\}$, $\rho(p') = b$ has a solution in **T**. Item 3 now follows from the fact that the mapping $\{p\} \to \{\rho(p'), b\}$ is linear in the size of instances.

For the backward direction of item 4 observe that a polynomial $p$ over $\overline{\mathbf{S}}$ is identically 0 if and only if for every $b \in \overline{S} \setminus \{0\}$, $p = b$ has no solution. To determine if $p$ is identically 0 requires that we refer to $|I||\Lambda| = mn$ instances of POL-SAT($\overline{\mathbf{S}}$); the size of each such instance is one more than the size of $\{p\}$. For the forward direction of Item 4, let $\{p, b\}$ be an instance of POL-SAT($\overline{\mathbf{S}}$). Let $b = [i, \lambda]$ ($i \in I, \lambda \in \Lambda$). Using Lemma 2.1, item 4, one can see that $p = b$ has a solution if and only if there exists a partial evaluation $e$ with domain $\{p_l, p_r\}$, such that $e(p_l) = [i, \alpha]$, $e(p_r) = [r, \lambda]$ ($\alpha \in \Lambda$, $r \in I$) and satisfying $e(p)$ is not identically 0. Thus to resolve the question of whether $p = b$ has a solution, it suffices to consider $|\Lambda||I| = mn$ instances of the POL=0 problem; the size of each such instance (the length of $e(p)$) is one less than the size of $\{p, b\}$. Item 4 now follows. □

Let **R** and **S** be Rees matrix semigroups. Let $p$ be a polynomial over **R**. If $\phi : \mathbf{R} \to \mathbf{S}$ be an isomorphism. Then the extension of $\phi$ which sends 1 of **R** to 1 of **S**, is a semigroup isomorphism also. Also $\phi$ induces a map of the polynomials over **R** ($\mathbf{R}^1$) to the polynomials of **S** ($\mathbf{S}^1$). To keep the notation brief we denote the map on polynomials also by $\phi$, where $\phi(p)$ maps each occurence of a constant symbol $s$ in $p$ to its image $\phi(s)$ and each variable occuring in $p$ is mapped to itself. Observe that if $p, q$ are polynomials and $b \in S$ ($b \in S^1$), then $p = q$ if and only if $\phi(p) = \phi(q)$ and $p = b$ has a solution if and only if $\phi(p) = \phi(b)$ has a solution. More generally, if $e$ is an evaluation of the variables $U$ occuring in $p$, letting $\phi(e) = \phi \circ e$, we have $e(p) = b$ if and only if $\phi(e)(\phi(p)) = \phi(b)$. These remarks also apply with $\mathbf{R}^1$ and $\mathbf{S}^1$ in place of **R** and **S** respectively.

Let $\mathbf{S} = \mathcal{M}(G; M)$ be a Rees matrix semigroup. Consider $M^T$, the transpose of $M$. The following notation is non-standard. We let $\mathbf{S^T}$ denote the Rees matrix semigroup $\mathcal{M}(G; M^T)$. For $[i, g, \lambda] \in S$, let $[i, g, \lambda]^T = [\lambda, g, i] \in S^T$. The map $s \to s^T$ induces a map of polynomials much as in the previous paragraph. For a polynomial $p$, let $p^T$ be the result of changing each occurence of a constant $s \in S$ to a constant $s^T \in S^T$. For two polynomials $p, q$ over **S** and $b \in S$, we have $p = q$ if and only if $p^T = q^T$ in $\mathbf{S^T}$ and $p = b$ has a solution if and only if $p^T = b^T$ has a solution in $\mathbf{S^T}$. Letting $1^T = 1$, the remarks of this paragraph hold with $\mathbf{R}^1$ and $\mathbf{S}^1$ in place of **R** and **S** respectively. The following lemma is an immediate consequence of the preceding two paragraphs.

**Lemma 2.6.** *If* **R** *and* **S** *are isomorphic Rees matrix semigroups or if* **R** *is isomorphic to* $\mathbf{S^T}$*, then if PR is one of the semigroup-complexity problem defined here, then PR(**R**)$\equiv_P$ PR(**S**). The same statement holds with* $\mathbf{R}^1$ *and* $\mathbf{S}^1$ *in place of* **R** *and* **S**, *respectively.*

# 3 Proofs for Theorems 1.3-1.5

We will intrepret the NP-complete GRAPH 3-COLORABILITY for simple graphs into semigroup-complexity problems for certain matrix semigroups. See [7] for background on the GRAPH 3-COLORABILITY problem. All graphs in this section are assumed to be simple and connected.

Let $G$ be a simple graph on $n$ vertices, $v_1, \ldots, v_n$. Let $\mathcal{P} = a_1, a_2, \ldots, a_p$ be a path in $G$ such that for $i$, $1 \leq i \leq p-1$, we let $a_i, a_{i+1}$ be an edge of $G$. We require that each edge of $G$ occurs in $\mathcal{P}$ in both directions. For any connected graph $G$ on $n$ vertices, there does in fact exist a path $\mathcal{P}$ of length $2n$ satisfying the condition described in the previous sentence. Let $t_G^{\mathcal{P}}(x_1, \ldots, x_n)$ $= t_G(x_1, \ldots, x_n) = \prod_i x_i$, where the order of the indices of the variables in the semigroup term $t_G$ is the same as the ordering of the indices of the vertices of the path $\mathcal{P}$.

Let $\mathbf{S} = \mathbf{S_{H_3}}$. In this paragraph we transform $t_G$ into a polynomial $\sigma(t_G)$ over $\mathbf{S_{H_3}}$, where $\sigma(t_G)$ involves added variables and constants. It will turn out that $\sigma(t_G)$ is identically equal to 0 over $\mathbf{S_{H_3}}$ if and only if $G$ can be 3-colored. The description of $\sigma(t_G)$ is in several steps. For each triple $(i, s, t)$, $1 \leq i \leq 3$, $1 \leq s \neq t \leq 3$, we introduce the polynomial $t_{(i)}^{(s,t)} = x_i (x_i^{(s,t)})^2 [s,t] x_i^{(s,t)} x_i$. Notice that $t_{(i)}^{(s,t)}$ involves the new variable $x_i^{(s,t)}$. For $1 \leq i \leq n$, let $\sigma(x_i)$ be the polynomial $x_i (\prod_{s,t} y_i^{(s,t)} t_{(i)}^{(s,t)} z_i^{(s,t)}) x_i$, where $y_i^{(s,t)}$ and $z_i^{(s,t)}$ are variables. We call $y_i^{(s,t)}$ and $z_i^{(s,t)}$ "buffer" variables. Their role will be to prevent the occurence of successive $x_i$'s in $\sigma(x_i)$. Observe that the number of variables and constants in $\sigma(x_i)$ is 50. Finally, we let $\sigma(t_G) = \prod_i \sigma(x_i)$. Note that as $G$ varies over the finite graphs, the size of the instance $\{\sigma(t_G)\}$ of POL=0($\mathbf{S_{H_3}}$) is linear in the size of the instance $\{G\}$ of GRAPH 3-COLORABILITY.

An element $s$ of a Rees matrix semigroup $\mathbf{S}$ is *nil* if $s^2 = 0$. Let $N$ be the set of nil elements of $\mathbf{S_{H_3}}$. Observe that by the definition of the multiplication operation of $\mathbf{S_{H_3}}$, we have $N = \{[z, z] : 1 \leq z \leq 3\}$.

**Lemma 3.1.** *$G$ is 3-colorable if and only if $\sigma(t_G)$ is not identically 0.*

*Proof.* We first describe the role of $t_{(i)}^{(s,t)}$ ($1 \leq s \neq t \leq 3$) in the construction. Let $f$ be a partial evaluation of the variables of $t_{(i)}^{(s,t)}$ with domain $\{x_i\}$. We claim that that if $f(x_i) = [t, s]$ and if $e$ is a total evaluation of the variables of $t_{(i)}^{(s,t)}$ extending $f$, then $e(t_{(i)}^{(s,t)}) = 0$. Conversely, if $f(x_i) \neq [t, s]$ and $f(x_i) \neq 0$, then there exists a total evaluation $e$ of the variables of $t_{(i)}^{(s,t)}$, extending $f$ and which satisfies $e(t_{(i)}^{(s,t)}) \neq 0$.

We prove the first part of the claim. Suppose that $e$ is an evaluation and $e(x_i) = [t, s]$. Assume for contradiction that $e(t_{(i)}^{(s,t)}) \neq 0$. Observe that the first coordinate of $e(x_i^{(s,t)})$ must be different from $s$ because $x_i x_i^{(s,t)}$ syntactically divides $t_{(i)}^{(s,t)}$; similarly, the second component must be different from $t$. Since $[s,t] x_i^{(s,t)}$ syntactically divides $t_{(i)}^{(s,t)}$, the first coordinate of $e(x_i^{(s,t)})$ cannot be $t$; similarly, the second component must be different from $s$. So $e(x_i^{(s,t)})$ is be of the form $[z,z]$ ($1 \leq z \leq 3$). But $(x_i^{(s,t)})^2$ syntactically divides $t_i^{(s,t)}$ and $[z,z] \in N$. So $e(t_{(i)}^{(s,t)}) = 0$ afterall. This completes the proof of the first part of the claim.

Suppose that $f(x_i) = [a,b] \neq [t,s]$ ($[a,b] \in S_{H_3} \setminus \{0\}$). There is enough symmetry in $t_{(i)}^{(s,t)}$ so that without loss of generality we may assume that $a \neq t$. We map $x_i^{(s,t)}$ to $[c,t]$ where $c \notin \{b,t\}$. In particular, $[c,t]^2 \neq 0$. We extend $f$ to $e$: $e(x_i) = [a,b]$ and $e(x_i^{(s,t)}) = [c,t]$. Observe that $e(t_{(i)}^{(s,t)}) \neq 0$.

We examine $\sigma(x_i)$. For any evaluation $e$ of the variables of $\sigma(x_i)$, if $e(x_i) \notin N$, then by the forward direction of the claim above, observe that $e(\sigma(x_i)) = 0$. If $f$ is a partial evaluation with domain $\{x_i\}$ and $f(x_i) \in N$, then by the converse direction of the previous claim, using the fact that the buffer variables separate possible repeated occurences of $x_i$, we can find a total evaluation $e$ (on the variables of $\sigma(x_i)$) which extends $f$ and satisfies $e(\sigma(x_i)) \neq 0$. Since $x_i$ is both the left-most and right-most symbol of $\sigma(x_i)$, by Lemma 2.1, item 4, $e(\sigma(x_i)) \neq 0$ implies $e(\sigma(x_i)) = [e(x_i)_1, e(x_i)_2] = e(x_i)$. For $\mathbf{S_{H_3}}$ we have shown the following.

**Observation 3.2.** *In the notation of the previous paragraph, $e(\sigma(x_i)) \neq 0$ implies $e(x_i) \in N = \{[z,z] : 1 \leq z \leq 3\}$ and $e(\sigma(x_i)) = e(x_i)$. Conversely, for any partial evaluation $f$ with domain $\{x_i\}$, if $f(x_i) \in N$, then there exists an evaluation $e$ of the variables of $\sigma(x_i)$, which extends $f$ and satisfies $e(\sigma(x_i)) = e(x_i)$.*

We examine $\sigma(t_G)$ and its connection with $G = (V, E)$. We are interested in whether $\sigma(t_G)$ is identically 0, so we assume for the remainder of the paragraph that $e$ is an evaluation of the variables of $\sigma(t_G)$ and that $e(\sigma(t_G)) \neq 0$. By Observation 3.2, we have that if $e(\sigma(t_G)) \neq 0$ and $x_i$ is a variable of $t_G$ (and thus a variable of $\sigma(t_G)$), then $e(\sigma(x_i)) = e(x_i) \in N$. Observe that $v_i, v_j$ is an edge of $G$ if and only if $x_i x_j$ syntactically divides $t_G$. Observe also that $x_i x_j$ syntactically divides $t_G$ implies $\sigma(x_i)\sigma(x_j)$ syntactically divides $\sigma(t_G)$. By the first sentence of Observation 3.2, we have $e(\sigma(t_G)) \neq 0$ implies $e(x_i) \neq e(x_j)$. In particular if $e(\sigma(t_G)) \neq 0$, then $e$ induces a proper 3-coloring of $G$, $c : V \to \{1,2,3\}$ as follows: for a vertex $v_i \in V$, $c(v_i) = z$ if $e(x_i) = [z,z]$ ($1 \leq z \leq 3$).

Conversely, if $c$ is a proper 3-coloring of $G$, $v_i$ is a vertex of $G$ and $c(v_i) = z$ ($1 \leq z \leq 3$), then we define an evaluation $f$ of the variables of $t_G$ in the obvious way: $f(x_i) = [z,z]$. By the second sentence of Observation 3.2, it follows that we can extend $f$ to $e$, a total evaluation of the variables of $\sigma(t_G)$ such that $e(\sigma(t_G)) \neq 0$. □

**Lemma 3.3.** *POL=0($\mathbf{S_{H_3}}$) is co-NP-complete.*

*Proof.* The size of $\{\sigma(t_G)\}$ in POL=0($\mathbf{S_{H_3}}$) is linear in the size of the instance $G$ of the 3-coloring problem. Apply Lemma 3.1. □

**Lemma 3.4.** *For $n \geq 4$, POL=0($\mathbf{S_{H_n}}$) is co-NP-complete.*

*Proof.* Suppose $n \geq 4$ and let $K = \{k : n \geq k \geq 4\}$. For the variable $x_i$, let $\alpha(x_i) = x_i(\prod_{u \in K}[u,u])x_i$. Observe that for any $[r,s]$ in $S_{H_n}$, $\alpha([r,s]) \neq 0$ if and only $1 \leq r \leq 3$ and $1 \leq s \leq 3$. Also $\alpha([r,s]) \neq 0$ implies that $\alpha([r,s]) = [r,s]$.

Let $p$ be a polynomial over $S_{H_3}$. Let $\alpha(p)$ be the polynomial over $S_{H_n}$ that results by substituting $\alpha(x_i)$ for each variable $x_i$ occuring in $p$. We have $p$ is identically 0 over $S_{H_3}$ if and only if $\alpha(p)$ is identically 0 over $S_{H_n}$. Since the mapping of the sizes of instances $\{p\} \to \{\alpha(p)\}$ is linear, the lemma follows. □

**Lemma 3.5.** *For any field $\mathcal{F}$, POL=0($L_2$ $\mathcal{F}$) is co-NP-complete.*

*Proof.* By Lemma 2.3, $(L_2\ \mathcal{F})/\mathcal{H}$ is isomorphic to $S_{H_n}$. Apply Lemma 2.5, item 1 (which states that the complexity of POL=0 does not change under moding out by $\mathcal{H}$), and Lemma 3.4. □

**Lemma 3.6.** *POL=0($L_n$ $\mathbf{Z_2}$) is co-NP complete for $n \geq 3$.*

*Proof.* Let $V = Z_2^n$. From the proof of Lemma 2.2, recall that $L_n\ \mathbf{Z_2}$ is isomorphic to a combinatorial Rees matrix semigroup. By the proof of Lemma 2.2, we can represent a non-0 element $A$ of $L_n\ \mathbf{Z_2}$ as a 2-tuple $[u,v]$, where $u,v \in V \setminus \{0\}$, $u$ spans the range of $A$, and the orthogonal complement of $v$ is the kernel of $A$. For two non-0 elements $[u,v], [u',v']$, we have $[u,v][u',v'] = [u,v']$ if $v^T u' = 1$ and $[u,v][u',v'] = 0$ otherwise.

Let $e_1, e_2, \ldots, e_n$ be a basis of $V$ and let $B$ be the set of non-zero vectors in subspace spanned by $\{e_1, e_2\}$. Let $x_i$ be a variable. For $a, b \in V$ satisfying $a+b \notin B$, consider the polynomial $q = q(x_i, y_i^{(a,b)}) = x_i y_i^{(a,b)}[a,a]y_i^{(a,b)}[b,b]y_i^{(a,b)}x_i$.

We claim that for an evaluation $e$ of the variables of $q$, if $e(x_i) = [a+b, v]$ or $e(x_i) = [v, a+b]$ ($v \in V$), then $e(q) = 0$. The claim follows from the observation that if $e(y_i^{(a,b)}) = [s,t]$ ($s,t \in V$) and $e(q) \neq 0$, then $1 = a^T s = a^T t = b^T s = b^T t$, from which it follows also that $(a+b)^T s = 0 = (a+b)^T t$. On the other hand, if $e(x_i) = [u,v]$ ($u,v \in V \setminus \{0\}$) and $a+b \notin \{u,v\}$, then it is not difficult to verify that there exists a choice for $e(y_i^{(a,b)})$ such that $e(q) \neq 0$. For each variable $x_i$ we let $\tau(x_i) = x_i(\prod_{a+b \notin B} u_i^{(a,b)} x_i y_i^{(a,b)}[a,a]y_i^{(a,b)}[b,b]y_i^{(a,b)}x_i v_i^{(a,b)})x_i$, where $u_i^{(a,b)}$ and $v_i^{(a,b)}$ are buffer variables. Let $p$ be an instance of POL=0($L_2$ $\mathbf{Z_2}$). We replace $x_i$ by $\tau(x_i)$ for each variable $x_i$ occuring in $p$ and denote the resulting polynomial by $\tau(p)$. Regard $\tau(p)$ as an instance of POL=0($L_n$ $\mathbf{Z_2}$). Observe that $p$ is identically 0 over $L_2\ \mathbf{Z_2}$ if and only $\tau(p)$ is identically 0 over $L_n\ \mathbf{Z_2}$. Again the mapping on size of instances is linear. The lemma follows from the fact that $L_2\ \mathcal{Z}_\in$ is isomorphic to $S_{H_3}$) and by Lemma 3.3. □

**Lemma 3.7.** *Let $\mathcal{F}$ be a field different from $\mathbf{Z_2}$. For $n \geq 2$, POL=0($L_n\mathcal{F}$) is co-NP complete.*

*Proof.* By Lemma 3.5 we can assume that $n \geq 3$. Let $V = \mathcal{F}^n$ and let $v_1, v_2, \ldots, v_n$ be an orthogonal basis of $V$. Let $C$ be a set of vectors in $V$ and let $R(C)$ denote the set of vectors in $V$ which are non-orthogonal to each of the vectors of $C$. Via a rather lengthy computation involving nothing more than dot products, the reader can verify that $\mathcal{F}$ is not $\mathbf{Z_2}$, then $R(R(v_1, v_2, v_1 + v_2))$ $= \langle v_1 \rangle \cup \langle v_2 \rangle \cup \langle v_1 + v_2 \rangle$, the union of the 1-dimensional subspaces of $V$ spanned by $v_1$, $v_2$, and $v_1 + v_2$.

We claim that we can choose an orthogonal basis so that for all $v \in \{v_1, v_2, v_1 + v_2\} = T$, we have that $v$ is orthogonal to exactly of one vector in $T$. To see this, observe that if $\mathcal{F}$ has characteristic 2, then the claim is trivial: use the standard ordered basis. Otherwise, use the fact that for any odd prime $p$, the equation $1 + c^2 + d^2 = 0 \bmod p$ has a solution. Let $e_1, \ldots, e_n$ denote the standard ordered basis. Let $v_1 = e_1$, $v_2 = ce_2 + de_3$, where $c, d$ satisfy $1 + c^2 + d^2 = 0 \bmod p$, and extend $\{v_1, v_2\}$ to an orthogonal basis of $V$. This completes the proof of the claim.

Note that to $T \times T$ can be associated a subsemigroup of $L_n \mathcal{F}$, namely the 0 matrix along with those rank 1 matrices $A$ such that the image of $A$ is a subspace spanned by a vector in $T$ and the kernel of $A$ is a subspace spanned by a vector in $T$. Call this subsemigroup $\mathbf{S}$. Of course $\mathbf{S}$ is isomorphic to a Rees matrix subsemigroup of $L_n \mathcal{F}$. We describe a structure matrix for $\mathbf{S}$. Let $v_1 = w_1$, $v_2 = w_2$ and $w_3 = v_1 + v_2$, let $N$ be the $3 \times 3$ matrix over $\mathcal{F}^*$ such that $N(r, \gamma) = w_\gamma^T w_r$ ($1 \leq \gamma \leq 3$, $1 \leq r \leq 3$). As we have mentioned, for all $v \in \{v_1, v_2, v_1 + v_2\} = T$, we have that $v$ is orthogonal to exactly one vector in $\{v_1, v_2, v_1 + v_2\}$. Hence $\mathbf{S}/\mathcal{H}$ is isomorphic to $\mathbf{S_{H_3}}$.

Now let $R = R(v_1, v_2, v_1+v_2))$ and let $\zeta(x_i) = x_i(\prod_{v \in R, w \in R} y_i^{(v,w)} x_i[v,w] x_i z_i^{(v,w)}) x_i$. For any evaluation $e$, we have if $e(x_i) \in T \times T$, then $e(\zeta(x_i)) = e(x_i)$; if $e(x_i)$ is not in $T \times T$, then $e(\zeta(x_i)) = 0$.

We have POL=0($\mathbf{S_{H_3}}$) $\equiv_P$ POL=0($\mathbf{S}$) $\leq_P$ POL=0($L_n \mathcal{F}$) $\leq_P$ POL-EQ($L_n \mathcal{F}$), the second to last inequality following from the previous paragraph and the last inequality a consequence of Lemma 2.5, item 2. The lemma follows from Lemma 3.3 which states that POL=0($\mathbf{S_{H_3}}$) is co-NP-complete. □

*Proof.* **Proof of Theorem 1.3** Theorems 1.3.1 follows from Lemma 3.6 and Lemma 3.7. From the proof of Lemma 3.7, we can infer also that POL=0($L_n \mathcal{F}$) is co-NP-complete. Now Theorem 1.3.2 follows from Lemma 2.5, item 2. □

*Proof.* **Proof of Theorem 1.4** It is easily verified that $L_n \mathcal{F}$ is a completely 0-minimal ideal of $T_n \mathcal{F}$. Theorem 1.4.1 now follows from an application of Lemma 2.5, item 2; Theorem 1.4.2 follows from an application of Lemma 2.5, item 3. □

As we have mentioned, Theorem 1.5 is a corollary of Theorem 1.4.

We prove Corollary 1.6. As we have mentioned, the POL-EQ$_\Sigma$ and POL-SAT$_\Sigma$ problems for any finite field are in P from which it follows easily that the same holds for any product of finite fields.

*Proof.* (Proof of Corollary 1.6) Let $\mathbf{R}$ be a finite non-commutative semi-simple ring; hence, $\mathbf{R}$ has a factor, say $\mathbf{R_i}$, which is a finite non-commutative simple ring. So $\mathbf{R_i}$ is isomorphic to $M_n \mathcal{F}$, where $n \geq 2$ and $\mathcal{F}$ is a field. We can assume that $\mathbf{R_i}$ is actually equal to $M_n \, \mathcal{F}$, where $n \geq 2$. To prove the corollary, it suffices, by Theorem 1.4, to polynomially intrepret POL=0($T_n \, \mathcal{F}$) into POL-EQ($\mathbf{R}$). Let $p$ be a polynomial over the semigroup $T_n \, \mathcal{F}$. Suppose that $v$ is a variable not occuring in $p$ and let $p' = cyp$ where $c$ is a non-0 element of $L_n \, \mathcal{F}$. Observe that $p$ is identically 0 over $T_n \, \mathcal{F}$ if and only if $p'$ is identically 0 over $T_n \, \mathcal{F}$. Of course $p'$ can be regarded as a polynomial (or a monomial) over $M_n \, \mathcal{F}$. We define $q$ a monomial over $\mathbf{R}$ by replacing each constant symbol $c$ occuring in $p'$ by the constant symbol $(0, \ldots, c, \ldots, 0) \in R$ where $c$ occurs in the i-th place. Since at least one constant occurs in $p'$ it follows that $p'$ is identically 0 over $T_n \, \mathcal{F}$ if and only if $q$ is identically equal to 0 over the ring $\mathbf{R}$. □

## 4 Proofs of Theorems 1.7, 1.8 and 1.10

In this section, all Rees matrix semigroups will be assumed to be combinatorial Rees matrix semigroups and $M$ will be refer to a regular 0-1 matrix. For $s = [i, \lambda]$ in Rees matrix semigroup $\mathbf{S_M}$, let $s_1 = i$ and $s_2 = \lambda$, the left-most and right-most components of $s$ respectively.

Our polynomial-time algorithms for polynomials and terms rely on a comparison of Z-sets of terms. We introduce another semigroup-complexity problem. We define TER=Z-SET (POL=Z-SET) for semigroups with a 0. An instance of TER=Z-SET (POL=Z-SET) is a pair of terms (polynomials) $\{p, q\}$ over $\mathbf{S}$. The question is whether $Z(p) = Z(q)$.

Let $\mathbf{S_M} = \mathbf{S}$ be a combinatorial Rees matrix semigroup with polynomial $p$ over $\mathbf{S}$. To $p$ we associate $B(p)$ a bipartite graph. As a first step in our description of $B(p)$, we define *bi-expressions*. Let $U$ be a set of variables. For $v \in U$, the bi-expression for $v$ is $(v_1, v_2)$. Let $U_1 = \{v_1 : v \in U\}$ and let $U_2 = \{v_2 : v \in U\}$. For $s = [i, \lambda] \in S$ ($i \in I$, $\lambda \in \Lambda$), the bi-expression for $s$ is $(i, \lambda) = (s_1, s_2)$. We use the directed graph $G(p)$ to define $B(p)$. Let $U$ be the variables involved in $p$ and let $K$ be the constants involved in $p$. We define a bipartite graph $B(p) = (X, Y, E)$ where $X = \{t_1 : t \in U \cup K\}$ and $Y = \{t_2 : t \in U \cup K\}$. For $i = 1$ or 2, we refer to $t_i$ as a *constant-vertex* if $t$ is a constant in $p$ and as a *variable-vertex* if $t$ is a variable. We define the edges of $B(p)$. For $v, w \in U \cup K$, $v_1, w_2$ is in $E$ if $v, w$ is an edge of $G(p)$. Of course $B(p)$ is bipartite since there are no edges among vertices of $X$ or among vertices of $Y$. We explain the relationship between $B(p)$ and $Z(p)$. Let $e$ be an evaluation of $U$ into $\mathbf{S}$. Each element of $U \cup S$ has a bi-expression. Observe that the evaluation $e$ determines maps $e_1 : U_1 \to I$ (for $v \in U$, $e_1(v_1) = e(v)_1$) and $e_2 : U_2 \to \Lambda$ (for $v \in U$, $e_2(v_2) = e(v)_2$). Let $\tilde{e} = e_1 \cup e_2$, $\tilde{e} : U_1 \cup U_2 \to I \cup \Lambda$, where $\tilde{e}|_{U_i} = e_i$, $i = 1, 2$.

**Observation 4.1.** *We have $e(p) = 0$ if and only if there exists variable $v$ occuring in $p$ and $e(v) = 0$, or there exist $x, y \in E$ ($E$ is the edge set for $B(p)$) such that $M(e_2(y), e_1(x)) = 0$ (equivalently, $M(\tilde{e}(y), \tilde{e}(x)) = 0$).*

A map $\tilde{e}$ where $\tilde{e} = e_1 \cup e_2$, $e_1 : U_1 \to I$ and $e_2 : U_2 \to \Lambda$, determines an evaluation $e : U \to S$ as follows. For $v \in U$, let $e(v) = [e_1(v_1), e_2(v_2)] = [\tilde{e}(v_1), \tilde{e}(v_2)]$.

The graphs $B(p)$ play a large role both here and in the companion paper ([14]) where connections between semigroup-complexity problems and certain complexity problems arising in a completely combinatorial setting are established. But we will need to define a modification of $B(p)$ in order to describe a polynomial time algorithm for POL=Z-SET($\mathbf{S_{I_k}}$), where $k \geq 2$ and $I_k$ is the $k \times k$ identity matrix.

Of course with $M = I_k$ we have $|I| = |\Lambda| = k$. It will be convenient to regard $I \cup \Lambda$ as the disjoint union of two copies of $I$. For a given polynomial $p$ over $\mathbf{S_{I_k}}$, if $[i, \lambda]$ occurs in $p$, then by the definition of $B(p)$ $i \in X$, $\lambda \in Y$ are vertices of $B(p)$. We denote elements of $X$ by $i_X$ and elements of $Y$ by $i_Y$ ($i \in I$). By $i$ we mean the pair $\{i_X, i_Y\}$ ($i \in I$). In this manner we identify the disjoint copies of $I$ and in so doing induce an identification of certain vertices of $B(p)$: if for some $i \in I$, $[r, i]$ and $[i, s]$ ($r, s \in I$) are constants occuring in $p$, then we identify the vertices $i_X \in X$ and $i_Y \in Y$. Denote the graph which results from under these identifications by $\overline{B}(p)$.

For example, let $u$ be a variable and consider the polynomial $p = [1, 1]u^2[1, 1]$ over $\mathbf{S_{I_2}}$. The reader can verify that $\overline{B}(p)$ is a triangle on the vertices $\{1, u_1, u_2\}$, where $(u_1, u_2)$ is the bi-expression for $u$. For any polynomial $p$, observe that in forming $\overline{B}(p)$ from $B(p)$, no variable-vertex of $B(p)$ is identified with any other vertex of $B(p)$. In particular, if $p$ is a term, since $B(p)$ has no constant-vertices, we have $\overline{B}(p) = B(p)$.

Denote the vertices of $\overline{B}(p)$ by $V$. Observe that for any evaluation $e$ of the variables of the polynomial $p$, we can compose $\tilde{e} : X \cup Y \to I \cup I$ (where $\tilde{e}$ is defined above, in connection with $B(p)$) with the map that identifies $i_X$ and $i_Y$ ($i \in I$) to form a map from $X \cup Y$ to $I$. Since $\tilde{e}(i_X) = \tilde{e}(i_Y)$ for any $i \in I$, we can go one step further and form a well-defined map $\hat{e} : V \to I$. Explicitly, if $v$ is a variable occuring in $p$ and $(e(v))_1 = i_X$, $(e(v))_2 = j_Y$ ($i, j \in I$), then $\hat{e}(v_1) = i$ and $\hat{e}(v_2) = j$. If $c = [i, j]$ ($i, j \in I$) is a constant occuring in $p$, then $\hat{e}(i_X) = i$ and $\hat{e}(j_Y) = j$. For example, for the polynomial $p = [1, 1]u^2[1, 1]$ of the previous paragraph, if $e(u) = [2, 1]$, then $\hat{e}(u_1) = 2$, $\hat{e}(u_2) = 1$, and $\hat{e}(1) = 1$. Note also that any map $\hat{f} : V \to I$ determines an evaluation $f$ of the variables of $p$ as follows. Let $u$ be a variable occuring in $p$ and let $f(u) = [\hat{f}(u_1), \hat{f}(u_2)]$.

Our analysis of POL=Z-SET for $\mathbf{S_{I_k}}$ will rely on comparisons of connected components of the graphs of the form $B(p)$ and $\overline{B}(p)$ where $p$ is a polynomial. We say that $C$ of $B(p)$ is *consistent* if $i_A$ and $j_B$ are vertices in $C$ ($i, j \in I$, $A, B \in \{X, Y\}$), then $i = j$. For example, letting $u, v, w$ be variables, for the polynomial $q$ over $\mathbf{S_{I_5}}$, $q = [1, 1]u^2[1, 2]v^2[3, 4]w[4, 4]w[5, 4]$, the components of $B(q)$ are $\{1_X, u_1, u_2, 1_Y\}$ ( a consistent component), $\{2_Y, v_1, v_2, 3_X\}$ ( not consistent), $\{4_Y, w_1\}$ (which is consistent), $\{w_2, 4_X, 5_X\}$ (which is not consistent).

**Observation 4.2.** *All components of $B(p)$ are consistent if and only if all for every component $D$ of of $\overline{B}(p)$, $D$ contains at most one element from $I$.*

**Lemma 4.3.** *For $k \geq 2$, POL=Z-SET($\mathbf{S_{I_k}}$) is in P. For polynomials $p, q$ over $\mathbf{S_{I_k}}$, to determine whether $Z(p) = Z(q)$ do the following. Compute the connected components of $B(p)$, $B(q)$, $\overline{B}(p)$, and $\overline{B}(q)$. We have $Z(p) = Z(q)$ if and only if one of the following hold:*

- $B(p)$ has a component which is not consistent and $B(q)$ has at a component which is not consistent.

- All components of $B(p)$ are consistent, all components of $B(q)$ are consistent and the components of $\overline{B}(p)$ and $\overline{B}(q)$ are the same.

Observe that for terms $p, q$ in the lemma above, we can replace all references to $\overline{B}$ with $B$. We demonstrate the need for $\overline{B}$ in analyzing the Z-sets of polynomials. Let $p = [1,1]u^2[1,1]$ and $q = [1,1]u[1,1]$. Observe that $Z(p) = Z(q)$ over $\mathbf{S_{I_2}}$. The components of both $B(p)$ and $B(q)$ are consistent, but $B(p)$ is a connected graph, while $B(q)$ is not. Note that $\overline{B}(p)$ and $\overline{B}(q)$ do in fact have the same components. In fact, they are both connected graphs over the same vertex set.

*Proof.* Let $p$ be a polynomial over $S_{I_k}$ and let $v, w \in U \cup K$, where $U \cup K$ denotes the symbols that occur in $p$. Recall that for an evaluation $e$ if $u$ is in $K$, then our convention is to let $e(u) = u$. We work first with $B(p)$. Suppose $C$ is a connected component of $B(p)$ and $v_1 \in X$ and $w_2 \in Y$ are contained in $C$. Then $M = I_k$ implies that for an evaluation $e$, if $e(p) \neq 0$, then $\tilde{e}(v_1) = \tilde{e}(w_2)$ ( use Fact 4.1 and induct on the distance between $v_1$ and $w_2$). In particular, if $C$ is a connected component of $B(p)$ and $C$ is not consistent, it follows that $p$ is identically 0.

We prove the converse: if all components of $B(p)$ are consistent, then $p$ is not identically 0. Since the components of $B(p)$ are consistent, we can construct a map $\hat{e}$ from the vertices of $\overline{B}(p)$ to $I$ so that $\hat{e}$ is locally constant, i.e. $\hat{e}$ is constant on the connected components of $\overline{B}(p)$. From $\hat{e}$ we construct an evalution $e$ of $U$. Let $v \in U$ be a variable. We let $e(v) = [\hat{e}(v_1), \hat{e}(v_2)]$. From Fact 4.1, it follows with a little thought that $e$ has been constructed so that $e(p) \neq 0$.

From a given polynomial $p$ over $\mathbf{S_{I_k}}$, we can construct $B(p)$ and determine if the connected components of $B(p)$ are consistent in polynomial time. Thus from the preceding two paragraphs we have shown the following.

**Fact 4.4.** *For $k \geq 2$, for a polynomial $p$ over $\mathbf{S_{I_k}}$, $p$ is identically 0 if and only if the components of $B(p)$ are consistent. In particular, POL=0($\mathbf{S_{I_k}}$) is in P.*

Moreover the preceding two paragraphs contain all the arguments that are needed for the following observation.

**Observation 4.5.** *For a polynomial $p$ over $\mathbf{S_{I_k}}$ and an evaluation $e$ of the variables $U$ of $p$, $e(p) \neq 0$ if and only if for each variable $v$ occuring in $p$, $e(v) \neq 0$ and for each component $D$ of $\overline{B}(p)$, $\hat{e}(D)$ is a singleton set.*

By Fact 4.4, if each of $B(p)$ and $B(q)$ has a component which is not consistent, then both $p$ and $q$ are identically 0; hence, $Z(p) = Z(q)$. If the components of $\overline{B}(p)$ and $\overline{B}(q)$ are the same, then from Observation 4.5, we have that that $Z(p) = Z(q)$. Thus if either of the two conditions described in the lemma holds, then $Z(p) = Z(q)$.

For the converse, assume that $Z(p) = Z(q)$. We prove that at least one of the two conditions described in the lemma holds. If both $p$ and $q$ are identically 0, then by Fact 4.4, each of $B(p)$ and $B(q)$ has a component which is not consistent. We may as well assume that the components of $B(p)$ and $B(q)$ are all consistent. Under this assumption, we claim that $Z(p) = Z(q)$ implies that the variables occuring in $p$ and $q$ are the same. Otherwise there exists an evaluation $e$ which send a variable $v$ occuring in $p$ but not $q$ to 0, while for each component $D$ of $\overline{B}(q)$, $\hat{e}(D)$ is a singleton set. By Observation 4.5, $e(q) \neq 0$.

Suppose for contradiction that there exist $a, b$, vertices in the same component of $\overline{B}(p)$, but $a, b$ are not in the same component of $\overline{B}(q)$. Let $D_a$, $D_b$ be the components of $a$ and $b$ respectively in $\overline{B}(q)$. Since all the components of $B(p)$ and $B(q)$ are consistent, we have that $D_a, D_b$ each contain at most one element from $I$. Since by assumption $D_a \neq D_b$, $D_a$ and $D_b$ contain no common element from $I$. Thus we can define a map $\hat{e}$ from the vertices of $\overline{B}(q)$ to $I$ so that $\hat{e}$ is constant on each connected component of $\overline{B}(q)$ and $\hat{e}(a) \neq \hat{e}(b)$. From $\hat{e}$, we define an evaluation $e$ of the variables involved in $p$ and $q$ in the usual way. Because $\tilde{e}$ is constant on each component of $\overline{B}(q)$, from Observation 4.5 we have $e(q) \neq 0$. Also from Observation 4.5, since $\hat{e}$ is not constant on components of $\overline{B}(p)$, it follows that $e(p) = 0$. We have shown that the algorithm described in the lemma is correct. Since the algorithm is clearly in P, we have POL=Z-SET($\mathbf{S_{I_k}}$) is in P. □

We have ignored the $\mathbf{S_{I_1}}$ case, since $\mathbf{S_{I_1}}$ is the 2 element semilattice. and it is easy to see that for semilattices all the semigroup-complexity problems defined here are in P. We will use the following fact.

**Fact 4.6.** *Two semilattice terms are equal if and only if they involve the same variables. For the two element semilattice $\mathbf{S}$, keeping in mind our assumption that 0 does not occur in any polynomial, two polynomials are equal if and only if they involve the same variables.*

Let $\mathbf{S}$ be a semigroup with 1. Since it will have no effect on the complexity of any semigroup-complexity problem associated with $\mathbf{S}$, for convenience we will assume that 1 does not occur in any polynomial over $\mathbf{S}$.

**Lemma 4.7.** *Let $k$ be a positive integer and let $p, q$ be polynomials over $\mathbf{S_{I_k}}^1$.*

1. *For $k = 1$, $Z(p) = Z(q)$ over $\mathbf{S_{I_k}}^1$ if and only if the variables involved in $p$ and $q$ are the same.*

2. *For $k \geq 2$, and terms $p, q$, $Z(p) = Z(q)$ over $\mathbf{S_{I_k}}^1$ if and only if at least one of the following holds.*

   (a) *Each of $B(p)$ and $B(q)$ has a component which is not consistent.*
   
   (b) *The components of $\overline{B}(p)$ and $\overline{B}(q)$ are the same.*

*In particular, POL=Z-SET($\mathbf{S_{I_k}}^1$) is in P.*

*Proof.* For the $k = 1$ case, we have $\mathbf{S_{I_k}}^1$ is the semilattice whose Hasse diagram is the the 3 element chain. For this case the lemma is trivial: for a polynomial $p$ and an evaluation $e$, $e(p) = 0$ if and only if there exists a variable $v$ occuring in $p$ such that $e(v) = 0$.

Assume $k \geq 2$ for the remainder of the proof. We extend the definition of bi-expression. Let $(u, u)$ be the bi-expression for the identity 1 of $\mathbf{S_{I_k}}^1$.

If $B(p)$ has a component which is not consistent, arguing as in the first paragraph of the proof of Lemma 4.3, we have $p$ is identically 0 over $\mathbf{S_{I_k}}^1$. Since $\mathbf{S_{I_k}}$ is a subsemigroup of $\mathbf{S_{I_k}}^1$, if the components of $B(p)$ are consistent, using Fact 4.4, $p$ is not identically 0. We have shown the following.

**Fact 4.8.** *For $k \geq 2$, for a polynomial $p$ over $\mathbf{S_{I_k}}^1$, we have $p$ is identically equal to 0 if and only if there exists a component of $B(p)$ which is not consistent. Thus, POL=0($\mathbf{S_{I_k}}^1$) is in P.*

Again using the fact that $\mathbf{S_{I_k}}$ is a subsemigroup of $\mathbf{S^1_{I_k}}$, it follows that for any pair of polynomials $p, q$ over $\mathbf{S^1_{I_k}}$, by Lemma 4.3, if $Z(p) = Z(q)$, then either the components of $\overline{B}(p)$ and the components of $\overline{B}(q)$ are the same, or each of $B(p)$ and $B(q)$ has a component which is not consistent.

We prove the converse. Assume that the components of $\overline{B}(p)$ and $\overline{B}(q)$ are the same. We claim that $Z(p) = Z(q)$ over $\mathbf{S_{I_k}}^1$. Let $e$ be an evaluation over $\mathbf{S_{I_k}}^1$. Modify $\overline{B}(p)$ and $\overline{B}(q)$ as follows. For each variable $v$ such that $e(v) = 1$, in each of $\overline{B}(p)$ and $\overline{B}(q)$ add the edge $v_1, v_2$. Call the resulting graphs $\overline{B}(p)_e$ and $\overline{B}(q)_e$ respectively. Since the components of $\overline{B}(p)$ and $\overline{B}(q)$ are the same, the components of $\overline{B}(p)_e$ and $\overline{B}(q)_e$ are the same. The following is a small modification of Fact 4.5 for $\mathbf{S_{I_k}}^1$: $e(p) \neq 0$ if and only if for each component $C$ of $\overline{B}(p)_e$, we have $\hat{e}(C) \setminus \{u\}$ has at most one element (we let $\hat{e}(u) = u$). Thus if the components of $\overline{B}(p)$ and $\overline{B}(q)$ are the same, then for an evaluation $e$ into $\mathbf{S_{I_k}}^1$ we have $e(p) = 0$ if and only if $e(q) = 0$.

In view of the preceding paragraphs, we have that the polynomial-time algorithm for POL=Z-SET($\mathbf{S_{I_k}}$) described in the previous lemma applies also to POL=Z-SET($\mathbf{S_{I_k}}^1$). □

Observe that the matrices $I_k$ are totally balanced. We extend Lemma 4.3 so that in modified form it can applied to any combinatorial Rees matrix semigroup with a totally balanced structure matrix. The following equivalent characterization of totally balanced matrices is useful. A matrix $M$ is totally balanced if and only if for all $i \neq j$ ($i, j \in I$) if there exists $\lambda \in \Lambda$, such that $M(\lambda, i) = 1 = M(\lambda, j)$, then $M^{(i)} = M^{(j)}$, and for all $\alpha \neq \beta$ ($\alpha, \beta \in \Lambda$), if there exists $r \in I$ such that $M(\alpha, r) = 1 = M(\beta, r)$, then $M_\alpha = M_\beta$.

Suppose a matrix $M$ has two identical rows (columns). The matrix $L$ that results by deleting one of the two equal rows (columns) will be said to result by *row retraction* (*column retraction*). We leave the proof of the following combinatorial lemma, and the observation that follows it, to the reader.

**Lemma 4.9.** *A regular $m \times n$ 0-1 matrix $M$ is totally balanced if and only if there exist a positive integer $k$ and a sequence of row and column retractions* applied to $M$ such that the matrix that results from the retractions is, up to row and column exchanges, equal to $I_k$, the $k \times k$ identity matrix.

In the preceding paragraph, we have $k = 1$ if and only if for some positive integers $m$ and $n$, $M$ is equal to the all 1's matrix $J_{m \times n}$.

**Observation 4.10.** *Let $M$ be a fixed regular totally balanced $m \times n$ matrix. The transformation of $p$, a polynomial over $\mathbf{S_M}$ ($\mathbf{S_M}^1$), into $\hat{p}$, the polynomial over $\mathbf{S_{I_k}}$ ($\mathbf{S_{I_k}}^1$) defined in the paragraph above, can be accomplished in linear time with respect to the length of $p$. To see why, note that the retraction aspect can be realized by two functions $t : \Lambda \to K_r$ and $s : I \to K_c$ and , where $K_r \subset \Lambda$, $K_c \subset I$ and $|K_r| = |K_c| = k$. The last aspect, involving mapping polynomials of a Rees matrix semigroup with a permutation matrix as structure matrix to one having a $k \times k$ identity matrix as structure matrix, is described by a pair of one-to-one maps $\phi_r$ and $\phi_c$, where $\phi_r : K_r \to K$ and $\phi_r : K_c \to K$, and $K = \{1, \ldots, k\}$. The maps $\phi_c \circ r$ and $\phi_r \circ s$ are not necessarily unique. Given a polynomial $p$ over $\mathbf{S_M}$, by replacing each constant $[i, \lambda]$ occuring in $p$ to $[\phi_c \circ r(i), \phi_r \circ s(\lambda)]$, the result is the polynomial $\hat{p}$.*

*In fact, given a regular totally balanced $m \times n$ matrix $M$, it is not hard to see that there exists on input $M$ a polynomial-time algorithm (with respect to $mn$) which outputs functions $\phi_c \circ r$ and $\phi_r \circ s$ that perform the transformation of $p$ into $\hat{p}$.*

Via Lemma 4.11 to follow, we will be able to extend Lemma 4.3 to semigroups of the form $\mathbf{S_M}$, where $M$ is totally balanced.

Let $M$ be a $m \times n$ 0-1 matrix. Suppose that that $M_\alpha$ and $M_\beta$ ($\alpha, \beta \in \Lambda$) are equal rows. Let $L$ be the $(m-1) \times n$ that results by deleting $M_\alpha$ from $M$. Let $p$ be a polynomial over $\mathbf{S_M}$. Let $p'$ be the polynomial which results by changing occurences of constants in $p$ of the form $[i, \alpha]$ to $[i, \beta]$ ($i \in I$). The variables occuring in $p$ and $p'$ are the same. Let $e$ be an evaluation of the variables of $p$. Let $e'$ be the evaluation defined as follows. If $v$ is a variable and $e(v) = [i, \alpha]$, let $e'(v) = [i, \beta]$. Otherwise, $e'(v) = e(v)$.

**Lemma 4.11.** *Suppose that $M$ is a regular 0-1 $m \times n$ matrix such that for $\alpha \neq \beta \in \Lambda$, $M_\alpha$ and $M_\beta$ are equal rows. Let $L$ be the matrix that results by deleting $M_\alpha$ from $M$. Then for polynomials $p, q$ over $\mathbf{S_M}$ ($\mathbf{S_M}^1$) we have $Z(p) = Z(q)$ if and only if $Z(p') = Z(q')$ over $\mathbf{S_L}$ ($\mathbf{S_L}^1$).*

*Statements of the same type hold for $M$ with column in place of row.*

*Proof.* Because $M_\alpha$ and $M_\beta$ are equal rows, it follows that for any polynomial $p$ over $\mathbf{S_M}$ ($\mathbf{S_M}^1$) and any evaluation $e$ of the variables of $p$, we have $e(p) = 0$ if and only if $e'(p') = 0$ over $\mathbf{S_L}$ ($\mathbf{S_L}^1$). As it pertains to rows of $M$, the lemma follows; use the second part of Lemma 2.6, a statement concerning transposes, to extend the result to the columns of $M$. □

Let $P$ be a $k \times k$ permutation matrix. Then by Lemma 2.1, item 5, there exists an isomorphism $\phi : \mathbf{S_P} \to \mathbf{S_{I_k}}$. As noted, $\phi$ induces a mapping between polynomials over $\mathbf{S_P}$ and $\mathbf{S_{I_k}}$.

Let $M$ be totally balanced. Since $M$ is totally balanced, by Lemma 4.9, there exists a sequence of row and column retractions which terminate in a permutation matrix $P$. Let $p$ be a polynomial over $\mathbf{S_M}$. We let $p^*$ be the polynomial over $\mathbf{S_P}$ that results from an application of a sequence of modifications of polynomials as described in the paragraph preceding Lemma 4.11. That is, as we retract $M$ sequentially at each stage we perform the "priming" operation on the the polynomial in question; the retraction procedure terminating in the permutation matrix $P$, with $p^*$ the end result of the associated sequence of primings. Since $P$ is a $k \times k$ permutation matrix, as noted above, there exists an isomorphism $\phi : \mathbf{S_P} \to \mathbf{I_k}$. Let $\hat{p} = \phi(p^*)$. Observe that for two polynomials $p, q$ over $\mathbf{S_P}$ ($\mathbf{S_P}^1$), we have $Z(p) = Z(q)$ if and only if $Z(\phi(p)) = Z(\phi(q))$ over $\mathbf{S_{I_k}}$ ($\mathbf{S_{I_k}}^1$). We use the notation of this paragraph in the next lemma.

**Lemma 4.12.** *Let $M$ be a totally balanced 0-1 matrix. Let $p, q$ be polynomials over $\mathbf{S_M}$ ($\mathbf{S_M}^1$).*

1. *If $M$ is an all 1's matrix, then $Z(p) = Z(q)$ over $\mathbf{S_M}$ ($\mathbf{S_M}^1$) if and only if the variables involved in $p$ and $q$ are the same.*

2. *If $M$ is not an all 1's matrix, then $Z(p) = Z(q)$ over $\mathbf{S_M}$ ($\mathbf{S_M}^1$) if and only if at least one of the following holds.*

   (a) *Each of $B(\hat{p})$ and $B(\hat{q})$ both have a component which is not consistent.*

   (b) *The components of $\overline{B}(\hat{p})$ and $\overline{B}(\hat{q})$ are the same.*

*In particular, for a totally balanced matrix $M$, POL=Z-SET($\mathbf{S_M}$) (POL=Z-SET(($\mathbf{S_M}^1$)) is in P.*

*Proof.* Item 1 is trivial since $M$ is a $J_{m \times n}$ matrix, for any polynomial $p$ over $\mathbf{S}$ ($\mathbf{S}^1$) and any evaluation $e$ of the variables of $p$, we have $e(p) = 0$ if and only if there exists a variable $v$ occuring in $p$ such that $e(v) = 0$.

By Lemma 4.11 we have $Z(p) = Z(q)$ over $\mathbf{S_M}$ ($\mathbf{S_M}^1$) if and only if $Z(p^*) = Z(q^*)$ over $\mathbf{S_P}$ ($\mathbf{S_P}^1$) where $P$ is a $k \times k$ permutation matrix, from which it then follows that $Z(p) = Z(q)$ over $\mathbf{S_M}$ ($\mathbf{S_M}^1$) if and only if $Z(\hat{p}) = Z(\hat{q})$ over $\mathbf{S_{I_k}}$ $\mathbf{S_{I_k}}^1$. Item 2 now follows by Lemma 4.3 (Lemma 4.7).

The last sentence of the lemma follows from items 1 and 2 of this lemma and the first paragraph of Observation 4.10. □

**Lemma 4.13.** *Let $M$ be a regular 0-1 matrix. If POL=Z-SET($\mathbf{S_M}$) is in P, then POL-EQ($\mathbf{S_M}$) and POL-SAT($\mathbf{S_M}$) are in P.*

*Proof.* The ideas are the much the same as in the proof of Lemma 2.5, item 4. Suppose $p, q$ are polynomials over $\mathbf{S_M}$. We have that $p = q$ if and only if $Z(p) = Z(q)$ and for all partial evaluations $f$ with domain $\{p_l, q_l, p_r, q_r\}$, if either $e(p_l)_1 \neq e(q_l)_1$ or $e(p_r)_2 \neq e(q_r)_2$, then $e(p)$ and $e(q)$ are both identically 0. To check if $p = q$ involves determining whether $Z(p) = Z(q)$ and then applying the POL=0 test described in Fact 4.4 on $2m^2n^2$ polynomials, each such instance of POL=0 of size the length of $p$ or the length of $q$. By Lemma 4.12, determining if $Z(p) = Z(q)$ is polynomial in the size of the instance $\{p, q\}$; by Fact 4.4, each of the $2m^2n^2$ tests can be be done in polynomial time with respect to the size of $\{p, q\}$. It follows that POL-EQ($\mathbf{S_M}$) is in P. Now from Lemma 2.5, item 4, we have POL-SAT($\mathbf{S_M}$) is also in P. □

*Proof.* **Proof of Theorem 1.10** Apply Lemma 4.12, followed by Lemma 4.13. □

Later for a totally balanced matrix $M$ we will provide explicit algorithms for the **term** equivalence problem for $\mathbf{S_M}$ and $\mathbf{S_M^1}$. First we consider the TER=Z-SET problem for not totally balanced combinatorial Rees matrix semigroups.

**Lemma 4.14.** *If $M$ is a regular 0-1, not totally balanced matrix and $p, q$ are* **terms** *over $\mathbf{S_M}$, then $Z(p) = Z(q)$ if and only if $G(p) = G(q)$. In particular, TER=Z-SET($\mathbf{S_M}$) is in P.*

*Proof.* Let $p, q$ be terms over $\mathbf{S_M}$. We claim that $G(p) = G(q)$ implies that $Z(p) = Z(q)$. Suppose that $G(p) \neq G(q)$. Let $x, y \in U$ be such that $x, y$ is an edge of $G(p)$ but $x, y$ is not an edge of $G(q)$. Since $M$ is not totally balanced, by Lemma 2.6, without loss of generality we can assume that $M(1, 1) = M(1, 2) = M(2, 1) = 1$ and $M(2, 2) = 0$. Let $e$ be the evaluation which sends $x$ to $[1, 2]$ and $y$ to $[2, 1]$ and which maps all other variables to $[1, 1]$. Observe that $x, y$ is an edge of $G(p)$ implies that $e(p) = 0$. Note that any product formed from $\{[1, 1], [2, 1], [1, 2]\}$ is 0 if and only if $[2, 1][1, 2]$ syntactically divide the product. Because $x, y$ is not an edge of $G(q)$, it follows that $e(q) \neq 0$.

We have proved that $Z(p) = Z(q)$ if and only if $G(p) = G(q)$. Computing $G(p)$ and $G(q)$ and determining if $G(p) = G(q)$ can be done in polynomial time. □

A specific polynomial-time algorithm for the term-equivalence problem for combinatorial Rees matrix semigroups is provided in the next lemma, one consequence of which is Theorem 1.7.

**Lemma 4.15.** *Let $\mathbf{S} = \mathbf{S_M}$ be a combinatorial Rees matrix semigroup and let $p, q$ be terms over $\mathbf{S}$.*

1. *If there exist positive integers $m, n$ such that $M = J_{m \times n}$ and $p, q$ are terms over $\mathbf{S_M}$, then $p = q$ if and only if*

   (a) *The variables occuring in $p$ and $q$ are the same.*
   
   (b) *If $m \geq 2$, $p_l = q_l$.*
   
   (c) *If $n \geq 2$, $p_r = q_r$.*

2. *If $M$ is totally balanced but $M \neq J_{m \times n}$ for any pair of positive integers $m, n$, then $p = q$ if and only if*

   (a) *The components of $B(p) = B(q)$ are the same.*

(b) $(p_l)_1$ and $(q_l)_1$ are the contained in the same connected component of $B(p) = B(q)$

(c) $(p_r)_2$ and $(q_r)_2$ are contained in the same connected component of $B(p) = B(q)$

(d) If there exist two equal rows of $M$, then $p_r = q_r$

(e) If there exist two equal columns of $M$, then $p_l = q_l$.

3. If $M$ is not totally balanced, then $p = q$ if and only if

   (a) $G(p) = G(q)$
   (b) $p_l = q_l$
   (c) $p_r = q_r$

*In particular, the term-equivalence problem for a combinatorial Rees matrix semigroup is in P. Moreover, on input $\{p, q, M\}$ where $p, q$ are both semigroup terms and $M$ is an $m \times n$ regular 0-1 matrix, we can determine if $p = q$ over $\mathbf{S_M}$ in polynomial time with respect to length(p)length(q)mn.*

*Proof.* Item 1 is trivial and follows from the fact that variables of $p$ and $q$ are the same if and only if $Z(p) = Z(q)$ and the definition of multiplication in the Rees matrix semigroup $\mathbf{S_M}$.

For item 2, we make use of the algorithm for POL=Z-SETS. For terms $p, q$ we have $Z(p) = Z(q)$ if and only if the components of $B(p)$ and the components of $B(q)$ are the same. Note that in the case of terms, the graphs $\overline{B}$ and $B$ are the same. Morever, since there are no occurences of constants in terms, components are automatically consistent.

If $(p_l)_1$ and $(q_l)_1$ are not contained in the same connected component, then there exists a map $\hat{e}$ which is constant on the components and which satisfies $\hat{e}((p_l)_1) \neq \hat{e}((q_l)_1)$. But $\hat{e}$ gives rise to an evaluation $e$ such that $e(p)_1 \neq e(q)_1$. An argument of the same type shows that $p = q$ implies that $(p_l)_1$, $(q_l)_1$ are contained in the same component.

If two rows of $M$ are equal, we can assume without loss of generality that $M_1 = M_2$. Suppose that $p_r \neq q_r$. By regularity of $M$, there exists $i \in I$ be such that $M(i, 1) = M(i, 2) = 1$. Let $e$ be an evaluation of the variables of $p, q$ satisfying $e(p_r) = [i, 1]$, $e(q_r) = [i, 2]$ and all other variables of $p, q$ are mapped to $[1, 1]$. Then $e(p) = [i, 1]$ while $e(q) = [i, 2]$. Thus $p = q$ implies that $e(p_r) = e(q_r)$. A similar argument shows that if two columns of $M$ are equal, then $p = q$ implies that $p_l = q_l$. This completes the forward direction of item 2.

Conversely suppose items 2 a-e hold for terms $p$ and $q$. By Lemma 4.12 Item 2a implies that $Z(p) = Z(q)$. For contradiction, assume that there exists an evaluation $e$ such that $e(p) \neq e(q)$. Since $Z(p) = Z(q)$, we have that $e(p) \neq 0$ and $e(q) \neq 0$. By Lemma 2.1, item 4, we must have that $e(p_l)_1 \neq e(q_l)_1)$ or $e(p_r)_r \neq e(q_r)_r$. Without loss of generality, assume that $e(p_l)_1 \neq e(q_l)_1$. From the preceding non-equality, we have $p_l \neq q_l$; hence, by item 2d the columns of $M$ must be distinct. By item 2b, $(p_l)_1$ and $(q_l)_1$ are contained in the same component. Suppose that $a_0 = (p_l)_1, \ldots, (q_l)_1 = a_s$ is a sequence witnessing that $(p_l)_1$ and $(q_l)_1$ are in the same component. Since $e(p) \neq 0$, we have $M(e_2(a_u), e_1(a_{u+1})) \neq 0$ ( $u$ is odd and $1 \leq u \leq s-1$). But $e(p_l)_1) \neq e(q_l)_1)$ implies that $s \geq 2$. In particular, there exists $\lambda \in \Lambda$, $i, j \in I$ such that $i \neq j$ and $M(\lambda, i) = 1 = M(\lambda, j)$. Since $M$ is totally balanced, it follows that $M^{(i)} = M^{(j)}$. But this contradicts the fact that no two columns of $M$ are the same. Thus $e(p_l)_1 = e(q_l)_1$ after all. A similar argument shows that $e(p_r)_2 = e(q_r)_2$. This completes the proof of item 2.

For item 3, assume $M$ is not totally balanced. Without loss of generality assume that $M(1,1) = M(1,2) = M(2,1) = 1$ and that $M(2,2) = 0$. Suppose that $p = q$. Thus $Z(p) = Z(q)$ and by by Lemma 4.14, we have $G(p) = G(q)$. Assume $x = p_l \neq q_l = y$. Consider the evaluation $e$ which sends $x$ to $[1,1]$, $y$ to $[2,1]$ and sends all other variables to $[1,1]$. Observe that $e(p) = [1,1]$ and $e(q) = [2,1]$. Thus, $p = q$ implies that $p_l = q_1$. By a similar argument, $p = q$ implies that $p_r = q_r$. For the converse, assume $G(p) = G(q)$, $p_l = q_1$, and $p_r = q_r$. By Lemma 4.14, we have $G(p) = G(q)$ implies $Z(p) = Z(q)$. By Lemma 2.1, item 4, for any evaluation $e$, if $e(p) \neq 0$, then $e(p) = [e(p_l)_1, e(p_r)_2] = [e(q_l)_1, e(q_r)_2] = e(q)$ (since $e(q) \neq 0$). This completes the proof of item 3.

The proof of the last paragraph is now easy. Given the matrix $M$, and semigroup terms $p$ and $q$, from $M$ one can decide in linear time (with respect to $mn$) which of the three polynomial-time algorithms described in the lemma should be applied to determine if $p = q$. □

We consider the TER=Z-SET problem for combinatorial Rees matrix semigroups with 1 adjoined.

**Lemma 4.16.** *If $M$ is a not totally balanced regular 0-1 matrix, then TER=Z-SETS($\mathbf{S_M}^1$) is in P.*

*Proof.* Let $p, q$ be terms over $\mathbf{S_M}^1$. If $p = q$, then because $p, q$ are terms, it is easy to see that the variables involved in $p$ and $q$ must be the same (otherwise send a variable $x$ occuring in $p$ but not $q$ to 0; send all the other variables to a non-0 idempotent). We assume the variables occuring in $p$ and $q$ are the same and denote them by $U$. Assume that $x, y \in U$. We define $B_{x,y,p}$, a collection of subsets of $U$. Let $U_i \subset U$ be in $B_{x,y,p}$ if there exists a polynomial $p_i$ such that $xp_iy$ syntactically divides $p$, neither $x$ nor $y$ occurs in $p_i$, and $U_i$ is the set of variables that occur in $p_i$.

Observe that $\emptyset$ is in $B_{x,y,p}$ if and only if $x, y$ is an edge of $G(p)$.

We define $A_{x,y,p}$ as follows. Under the partial ordering on the set $U$ determined by set inclusion, we let $A_{x,y,p}$ be the minimal sets in $B_{x,y,p}$. Observe that $A_{x,y,p} = \emptyset$ if and only if $x, y$ is an edge of $B_{x,y,p}$ and that $A_{x,y,p}$ is an anti-chain.

**Claim** *If $Z(p) = Z(q)$, then for all $x, y \in U$ $A_{x,y,p} = A_{x,y,q}$.*

If $A_{x,y,p} \neq A_{x,y,q}$, then without loss of generality there exists $D \in A_{x,y,p}$ such that for all $E \in A_{x,y,q}$, $D$ is not a subset of $E$. We may assume that $M(1,1) = M(1,2) = M(2,1) = 1$ and $M(2,2) = 0$. Consider the following evaluation $e$. Let $e(D) = \{1\}$. Let $e(x) = [1,2]$, $e(y) = (2,1)$ and let $e$ send the remaining variables to $[1,1]$. Observe that $e(p) = 0$. A product composed of elements in $\{[1,1], [1,2], [2,1], 1\}$ is 0 if and only if $[1,2]$ occurs to the left of $[2,1]$, separated by a sequence of 1's. By the first sentence of this paragraph, no such product can be realized in $e(q)$. Thus $e(q) \neq 0$.

**Claim** If for all $x, y \in U$ $A_{x,y,p} = A_{x,y,q}$, then $Z(p) = Z(q)$.

Assume for contradiction that for all $x, y \in U$ $A_{x,y,p} = A_{x,y,q}$ but $Z(p) \neq Z(q)$. Without loss of generality there exists an evaluation $e : U \to S_M^1$ such that $e(p) \neq 0$ and $e(q) = 0$. The hypothesis of the claim guarantees that $G(p) = G(q)$; hence, $e(p) = 0$ does not result in the following manner: $u, v$ is an edge of $G(p)$, $e(u), e(v) \in S_M$ and $e(u)e(v) = 0$. Thus $e(p) = 0$ results as follows: $uu_1 \ldots u_k v$ syntactically divides $p$ (where $u, v, u_1, \ldots, u_k \in U$), $e(u)e(v) = 0$ (so $e(u), e(v) \in S_M$), $e(u_i) = 1$ ($i = 1, \ldots, k$). By definition, we have $\{u_1, \ldots, u_k\} \in B_{u,v,p}$. By the hypothesis of this claim, there exists a subset of $\{u_1, \ldots, u_k\}$ in $A_{u,v,q}$. But then $e(q) = 0$, a contradiction. This completes the proof of the claim.

By the preceding claims, it follows that to check whether $Z(p) = Z(q)$ it suffices to verify that for all $x, y \in U$, $A_{x,y,p} = A_{x,y,q}$. This can be checked by a polynomial-time algorithm. □

The following two lemmas constitute the proof of Theorem 1.8.

Suppose that $p$ is a term over $\mathbf{S_M}^1$ and $v$ occurs in $p$. Let $p/v$ denote the term that results from the elimination of all occurences of $v$ in $p$. If $v$ occurs in two terms $p$ and $q$, then $p = q$ over $\mathbf{S_M}^1$ implies that $p/v = q/v$ over $\mathbf{S_M}^1$.

By the *left (right) p-sequencing* of the variables occuring in $p$, we mean the ordered set $\{x_{i_1}, \ldots, x_{i_n}\}$, which records the order in which the variables of $p$ make their appearance from left to right (right to left) in $p$. If $M$ is not totally balanced, then by the last sentence of the previous paragraph and Lemma 4.15, item 3, along with the fact that $\mathbf{S_M}$ is a subsemigroup of $\mathbf{S_M}^1$, it follows that if $p = q$ over $\mathbf{S_M}^1$, then the left (right) $p$-sequencing and the left (right) $q$-sequencing are the same.

**Lemma 4.17.** *Let $M$ be a not totally balanced 0-1 matrix and let $p, q$ be terms over $\mathbf{S_M}^1$. Then $p = q$ if and only if*

1. *The left (right) p-sequencing and q-sequencing are the same.*

2. *For all variables $x, y$ occuring in $p, q$, $A_{x,y,p} = A_{x,y,q}$.*

*In particular, if $M$ is not totally balanced, then TER-EQ($\mathbf{S_M}^1$) is in P.*

*Proof.* The forward direction follows from the remarks of the paragraph immediately preceding this lemma and the proof of Lemma 4.16. For the converse, suppose that items 1 and 2 above hold for terms $p$ and $q$. Since item 2 holds, by the proof of Lemma 4.16, we have $Z(p) = Z(q)$. Assume there exists an evaluation $e$ such that $e(p) \neq e(q)$. Observe that case $e(p) \neq 1$ and $e(q) \neq 1$. We have $p_l = q_l$ and $p_r = q_r$ by item 2. Suppose $e(p_l) \neq 1$ and $e(p_r) \neq 1$. Then $Z(p) = Z(q)$ and Lemma 2.1, item 4, imply that $e(p) = e(q)$. Hence we can assume $e(p_l) = 1$. Let $p'$ and $q'$ be the terms that result by eliminating $p_1$ from $p$ and $q$ respectively. Observe that $p'$, $q'$ satisfy items 1 and 2. By induction on length of instance, we can assume that $p' = q'$. But $e(p) = e(p') = e(q') = e(q)$ , contradicting that $e(p) \neq e(q)$.

We have shown that the algorithm above for term equivalence of the lemma is correct. It is easy to verify that the algorithm can be implemented in polynomial time.

□

Let $p$ be a semigroup term. The *left (right) component sequencing* of $p$ is the order in which representatives of the components of $B(p)$ make their appearance from left to right (right to left) in $p$.

**Lemma 4.18.** *Let $M$ be totally balanced. Let $p, q$ be terms over $\mathbf{S_M}^1$. If $M = J_{m,n}$, an all 1's matrix, then $p = q$ if and only if*

1. *The variables occuring in $p$ and $q$ are the same.*

2. *If $m \geq 2$, then the right p-sequencing and the right q-sequencing is the same.*

3. *If $n \geq 2$, then the left p-sequencing and the left q-sequencing is the same.*

*If $M$ is not an all 1's matrix, then $p = q$ if and only if*

1. *Each of $B(p)$ and $B(q)$ have components which are not consistent, or the components of $B(p)$ and $B(q)$ are the same.*

2. *The left component p-sequencing and and the left component q-sequencing are the same; the right component p-sequencing and the right component q-sequencing are the same.*

3. *If $M$ has two equal columns, then the left p-sequencing and the left q-sequencing are the same.*

4. *If $M$ has two equal rows, then the right p-sequencing and the right q-sequencing are the same.*

*In particular, if $M$ is totally balanced, then TER-EQ($\mathbf{S_M}^1$) is in P.*

*Proof.* Using the fact that $\mathbf{S_M}$ is a subsemigroup of $\mathbf{S_M}^1$, the forward direction follows from Lemma 4.15 and the remark following Lemma 4.16 (involving the elimination of a variable occuring in $p$ and $q$). For the converse, assume that items pertaining to the given case (all 1's matrix or otherwise) but that $p \neq q$. Let $e$ be an evaluation such that $e(p) \neq e(q)$. We claim that we may assume that $e(p_1) \neq 1$. If not we could eliminate all instances of the variable $p_1$ from $p$ and $q$. Denote the resulting terms $p', q'$. Observe that $p', q'$ still satisfy the hypotheses. By inducting on the number of variables, we may assume $p' = q'$ from which it follows that $e(p) = e(q)$ after all. By an argument of the same type, we may assume that $e(p_l) \neq 1$. Item 1 in each of the two cases implies that $Z(p) = Z(q)$ (using Lemma 4.12). Thus $e(p) \neq e(q)$ implies $e(p) \neq 0$ and $e(q) \neq 0$. In the all 1's case, items 2 and 3 guarantee that $\tilde{e}(p_l) = \tilde{e}(q_1)$ and that $\tilde{e}(p_r) = \tilde{e}(q_r)$. Now $e(p) = e(q)$ follows from Lemma 2.1, item 4, a contradiction. So $p = q$ for this case. In the non-all 1's case, with a little thought it can be seen that items 2-4 guarantee that $e(p) \neq 0$ and $e(q) \neq 0$ imply that $\hat{e}(p_l) = \hat{e}(q_l)$ and $\hat{e}(p_r) = \hat{e}(q_r)$, from which it follows that $e(p) = e(q)$ after all. □

## 5 Conclusion

In this section we hope to provide some impetus for further work on the term-equivalence and polynomial-equivalence problems for general algebras, particularly semigroups.

For a finite group **G** with exponent $n$, observe that that for all $x \in G$, $x^{n-1} = x^{-1}$. Hence we can replace any term $p$ over **G** by a term $p'$ such that $p'$ involves no instances of the inverse operation. Let $\mathcal{M}(G; M)$ be a Rees matrix semigroup and let $p$ be a semigroup term. By $p^G$ we mean $p$ intrepreted as a term over **G**.

Let $\mathbf{S} = \mathcal{M}(G; M)$, where $M$ is an $m \times n$ regular matrix over **G**, be a Rees matrix semigroup. Let $M(\lambda, i) = h \neq 0$ ($i \in I$, $\lambda \in \Lambda$). We will need some additional information concerning Rees matrix semigroups.

It is straight-forward to verify that $\{[i, g, \lambda] : g \in G\}$ is a subgroup of **S** which is isomorphic to **G** (under the map $[i, g, \lambda] \to gh$). Also, Lemma 2.1, item 6 states that if a row of $M$ is multiplied on the left (or a column on the right) by an element $g \in G$ resulting in a matrix we will call $L$, then $\mathcal{M}(G; M)$ is isomorphic to $\mathcal{M}(G; L)$. We will need to know the explicit isomorphism. The map $[i, h, \lambda] \to [i, k, \lambda]$, where $k = hg$, if $\lambda = \alpha$, and $k = h$ otherwise, describes an isomorphism $\rho: \mathcal{M}(G; M) \to \mathcal{M}(G; L)$. On the other hand if $r \in I$, $g \in G$, and we let $L$ be the $m \times n$ matrix which results by replacing $M^{(r)}$ by $M^{(r)}g$, then the map $[i, h, \lambda] \to [i, k, \lambda]$ where $k = gh$ if $i = r$ and $k = h$ otherwise defines an isomorphism $\rho: \mathcal{M}(G; M) \to \mathcal{M}(G; L)$. In particular, if $L$ results from any series of row and column "re-scalings", then there exists an isomorphism $\rho: \mathcal{M}(G; M) \to \mathcal{M}(G; L)$ such that $\rho$ carries $[i, h, \lambda] \to [i, k, \lambda]$. In particular, for any isomorphism $\rho$ associated with a sequence of column and row, we have $\rho$ acts identically on $I$ and $\Lambda$.

**Lemma 5.1.** *Let* $\mathbf{S} = \mathcal{M}(G; M)$, *where* $M$ *is a matrix over* **G**. *Then TER-EQ(***G***)$\leq_P$ TER-EQ(***S***).*

*Proof.* Let $p, q$ be terms over **S** with associated group terms $p^G, q^G$. We are looking for an algorithm that determines whether $p^G = q^G$ by somehow referring the question to TER-EQ(**S**). We will assume that variables in $p, q$ are the same: otherwise, if $y$ occurs in $p$ but not in $q$, multiply $q$ on the left by $y^t$, where $t$ is the exponent of **G**.

Choose a variable $x$ that does not occur in $p$. For each variable $x_i$ occuring in $p$, let $\zeta(x_i) = x^t x_i x^t$. Let $\zeta(p), \zeta(q)$ be the terms that result by replacing $x_i$ by $\zeta(x_i)$ for every variable $x_i$ occuring in $p$. Observe that $p^G = q^G$ if and only if $\zeta(p)^G = \zeta(q)^G$. Observe that for the graphs $G(\zeta(p)), G(\zeta(q))$, we have $G(\zeta(p)) = G(\zeta(q))$. Also, $\zeta(p)_l = x = \zeta(q)_l$ and $\zeta(p)_r = x = \zeta(q)_r$. By Lemma 4.15, which characterizes equal terms over a combinatorial Rees matrix semigroup, we have $\overline{\zeta(p)} = \overline{\zeta(q)}$. We claim that $\zeta(p)^G = \zeta(q)^G$ over $\mathbf{G}$ if and only $\zeta(p) = \zeta(q)$ over $\mathbf{S}$. One direction is trivial since $\mathbf{G}$ is isomorphic to a subgroup of $\mathbf{S}$. Conversely, assume for contradiction that $\zeta(p)^G = \zeta(q)^G$ but that $\zeta(p) \neq \zeta(q)$ over $\mathbf{S}$. As we have remarked, $\overline{\zeta(p)} = \overline{\zeta(q)}$; hence, $Z(\zeta(p)) = Z(\zeta(q))$. Thus, $\zeta(p) \neq \zeta(q)$ implies that there exists an evaluation $e$ of the variables of $\zeta(p)$ such that $e(\zeta(p)) \neq 0$, $e(\zeta(q)) \neq 0$ and $e(\zeta(p)) \neq e(\zeta(q))$. Let $e(x) = [i, g, \alpha]$. If the variable $v$ occurs in $p, q$, because we have jacketed each occurence of $v$ with occurences to the left and right of $x$, only the ith column of $M$ and the $\alpha$-row of $M$ are involved in computing $e(\zeta(p))$ and $e(\zeta(q))$.

We can perform row and column re-scalings of $M$ so that for the resulting matrix $L$ we have $L^i$ and $L_\alpha$ both consists of 0's and 1's. Associated with these rescalings is an isomorphism of Rees matrix semigroups, $\mu : \mathcal{M}(G; M) \to \mathcal{M}(G; L)$, as described in the paragraph preceding this lemma. Applying $\mu$ to the inequality $e(\zeta(p)) \neq e(\zeta(q))$ results in an inequality over the Rees matrix semigroup $\mathcal{M}(G; L)$ which we denote by $\mu(e)(\zeta(p)) \neq \mu(e)(\zeta(q))$. By the paragraph preceding the lemma, the left and right components of the 3-tuples of $\mathbf{S}$ remain invariant under the isomorphism $\mu$. Thus without loss of generality we may assume that there exists an evaluation $e$ such that $e(x) = [i, g, \alpha]$, $e(\zeta(p)) \neq 0$, $e(\zeta(q)) \neq 0$, $e(\zeta(p)) \neq e(\zeta(q))$, and $M_\alpha$ and $M^{(i)}$ both consist of 0's and 1's. The fact that $M_\alpha$ and $M^{(i)}$ both consist of 0's and 1's allows us to do the computations readily. If $v$ is a variable and $e(v) = [i, g, \lambda]$, let $e^G(v) = g$ and let $e^G(p^G)$ and $e^G(q^G)$ be the evaluations of the group terms $p^G$ and $q^G$ that result from substitutions of $v$ by $e^G(v)$. It is easy to see that $e(\zeta(p)) = [i, e^G(p^G), \alpha]$ and $e(\zeta(q)) = [i, e^G(q^G), \alpha]$. By assumption, $p^G = q^G$; hence, $e(\zeta(p)) = e(\zeta(q))$, a contradiction. Since the size of $\{\zeta(p), \zeta(q)\}$ is linear in the size of $\{p^G, q^G\}$, the lemma follows from the claim. □

**Lemma 5.2.** *Let* $\mathbf{S} = \mathcal{M}(G; M)$, *where* $M$ *is a 0-1 matrix. Then TER-EQ($\mathbf{G}$)$\equiv_P$ TER-EQ($\mathbf{S}$).*

*Proof.* Let $p, q$ be terms over $\mathbf{S}$. We claim that $p = q$ if and only $\overline{p} = \overline{q}$ and $p^G = q^G$. One direction of the claim is obvious. For the other direction, suppose that $\overline{p} = \overline{q}$ and $p^G = q^G$. By Lemma 4.15 $\overline{p} = \overline{q}$ implies that the variables occuring in $p$ and $q$ are the same and that $p_l = q_l$ and $p_r = q_r$. Suppose that $e$ is an evaluation of the variables of $p, q$ such that $e(p) \neq 0$. Note that $\overline{p} = \overline{q}$ implies that $e(q) \neq 0$. That the non-0 entries of the matrix $M$ are 1's implies that $e(p) = [(e(p_l)_1, e^G(p), (e(p_r)_2]$. For the same reasons, $e(q) = [(e(q_l)_1, e^G(q), (e(q_r)_2]$. The claim follows from $p_l = q_l$ and $p_r = q_r$. By Theorem 1.6, $\overline{p} = \overline{q}$ can be determined in polynomial time. From the claim above, it now follows that TER-EQ($\mathbf{S}$) $\leq_P$ TER-EQ($\mathbf{G}$), and the lemma now follows by an application of Lemma 5.1. □

The next problems are motivated by curiousity concerning the extent to which the term-equivalence problem for a semigroup $\mathbf{S}$ is determined by the complexities of the term equivalence for the subgroups of $\mathbf{S}$.

**Problem 5.3.** *For a Rees matrix semigroup* $\mathbf{S} = \mathcal{M}(G; M)$, *does TER-EQ(***G***) in P imply TER-EQ(***S***) is in P? Is it true that TER-EQ(***S***)$\equiv_P$ TER-EQ(***G***) for an arbitrary Rees matrix semigroup ***S***?*

A semigroup is said to be a *combinatorial* semigroup if all of its subgroups are trivial.

**Problem 5.4.** *Let* ***S*** *be a combinatorial semigroup. Is TER-EQ(***S***) in P?*

As mentioned at the end of Section 1.3.1, we work out TER-EQ($\Omega(\mathbf{S})$) for certain non-combinatorial Rees matrix semigroups ***S***. For the next lemma, familiarity with semigroup theory, particularly translational hull theory, is assumed. The lemma indicates that the term equivalence problem for translational hulls of Rees matrix semigroups with 0-divisors may be of some interest.

**Lemma 5.5.** *Let $M$ be the $2 \times 2$ matrix such that $M(1,1) = 0$, $M(i,j) = 1$ for $i \neq 1$ or $j \neq 1$, let ***G*** be a group and let $\mathbf{S} = \mathcal{M}(M; G)$. Then TER-EQ($\Omega(\mathbf{S})$) $\equiv_P$ TER-EQ(***G***).*

*Proof.* A computation using Petrich's description of the translational hull (See [13]) shows that $\Omega(\mathbf{S})$ consists of ***S*** along with the group of invertible elements of $\Omega(\mathbf{S})$, a group which is isomorphic to ***G***. Letting $g$ denote an invertible element of $\Omega(\mathbf{S})$. Then it turns out that from the Petrich description we can infer that for any $[i, h, \lambda] \in S$, we have $g[i, h, \lambda] = [i, gh, \lambda]$ and $[i, h, \lambda]g = [i, hg, \lambda]$ (so the invertible elements act identically on $I$ and $\Lambda$, a help in analyzing $\Omega(\mathbf{S})$). Observe that the Green's Relation $\mathcal{H}$ is a congruence of $\Omega(\mathbf{S})$ and that $\Omega(\mathbf{S})/\mathcal{H}$ is isomorphic to a combinatorial Rees semigroup (with a structure matrix which is not totally balanced) with 1 adjoined. The term-equivalence problem for such a semigroup is in P by Lemma 4.17. For a term $p$ over $\Omega(\mathbf{S})$, let $p/\mathcal{H}$ be denoted by $\overline{p}$.

We claim that for any pair of terms $p, q$ over $\Omega(\mathbf{S})$ that $p = q$ if and only if $\overline{p} = \overline{q}$ and $p^G = q^G$. Clearly $p = q$ implies $\overline{p} = \overline{q}$ and $p^G = q^G$. Conversely, assume $\overline{p} = \overline{q}$ and $p^G = q^G$. Let $e$ be an evaluation of the variables of $p$ and assume that $e(p) \neq 0$ (so $e(q) \neq 0$ also). For the evaluation $e$, we define $e^G$ in the context of $\Omega(\mathbf{S})$ as follows. If $e(v) = [i, g, \lambda]$, then $e^G(v) = g$ (just as in Lemma 5.1). If $e(v) = h$, an invertible element of $\Omega(\mathbf{S})$, let $e^G(v) = h$. Thus $e^G$ defines an evaluation of $p^G$ and $q^G$. Since the invertible elements of $\Omega(\mathbf{S})$ act identically on $I$ and $\Lambda$ and $M$ is a 0-1 matrix, it follows that $e(p) = [\overline{e(p)}_1, e^G(p^G), \overline{e(p)}_2]$ and $e(q) = [\overline{e(q)}_1, e^G(q^G), \overline{e(q)}_2]$. The claim follows since $\overline{p} = \overline{q}$ implies that $\overline{e(p)}_1 = \overline{e(q)}_1$ and $\overline{e(p)}_2 = \overline{e(q)}_2$ and that $p^G = q^G$ implies $e^G(p^G) = e^G(q^G)$.

From the claim and the fact that the equality of $\overline{p}$ and $\overline{q}$ can be determined in polynomial-time, it follows that TER-EQ($\Omega(\mathbf{S})$) $\leq_P$ TER-EQ(***G***). We reverse the preceding inequality. Let $p^G$ and $q^G$ be group terms. We may assume that the variables occuring in $p$ and $q$ are the same. Let $U = \{x_1, \ldots, x_r\}$ be the variables that occur in $p, q$. Let $t$ be the exponent of ***G***. Let $\sigma(x_i) = (x_1^t \ldots x_r^t) x_i (x_1^t \ldots x_r^t)$ ($1 \leq i \leq r$). Let the terms that results by replacing $x_i$ by $\sigma(x_i)$ for each variable $x_i$ occuring in $p$ and $q$ be denoted $\sigma(p)$ and $\sigma(q)$ respectively. Observe that for any pair of variables $x, y \in U$, $A_{x,y,\sigma(p)} = A_{x,y,\sigma(q)}$.

Moreover, a quick perusal of the statement of Lemma 4.17 should convince the reader that $\overline{\sigma(p)} = \overline{\sigma(q)}$. Since $t$ is the exponent of $\mathbf{G}$, it follows that $p^G = \sigma(p)^G$ and $q^G = \sigma(q)^G$. Applying the claim of the previous paragraph, we have $\sigma(p) = \sigma(q)$ over $\Omega(\mathbf{S})$ if and only if $\sigma(p)^G = \sigma(q)^G$ over $\mathbf{G}$ if and only if $p^G = q^G$ over $\mathbf{G}$. Since the mapping $\{p^G, q^G\} \to \{\sigma(p), \sigma(q)\}$ is polynomial in size, we have TER-EQ($\mathbf{G}$)$\leq_P$ TER-EQ($\Omega(\mathbf{S})$). □

At the end of the introduction to the paper we mentioned that there is a strong connection between problems for combinatorial Rees matrix semigroups and the graph retract problem for bipartite graphs. Let $M$ be a regular 0-1 $m \times n$ matrix. Then $M$ determines a bipartite graph $B_M = (X, Y, E)$ where the $X$ is the set of rows of $M$ (indexed as usual by $\Lambda$), $Y$ is the set of columns of $M$ (indexed by $I$) and for $i \in I$ and $\lambda \in \Lambda$, $\lambda, i$ is in $E$ if $M(\lambda, i) \neq 0$. In a sequel to this paper, we show that POL=0($\mathbf{S_M}$) is in P if and only if RET $B_M$ is in P, and that POL=0($\mathbf{S_M}$) is co-NP-complete if and only if RET $B_M$ is NP-complete. It is thought that the retraction problem for bipartite graphs may not admit a "dichotomy". That is, if P$\neq$ NP, then it may be that there exist bipartite graphs $B$ such that RET $B$ is neither in P nor is NP-complete (See [6]). Since POL=0 polynomially embeds into POL-EQ for combinatorial Rees matrix semigroups (See Lemma 2.5, item 2), it is possible that POL-EQ for combinatorial Rees matrix semigroups may not admit a dichtotomy.

Note that for a given algebra $\mathbf{A}$, by modifying $\mathbf{A}$ by adding to the constant functions to the fundamental operations of $\mathbf{A}$ and denoting the new algebra by $\mathbf{A}_A$, we have TER-EQ($\mathbf{A}_A$) and POL-EQ($\mathbf{A}_A$) have the same complexity. The lack of certainty in the preceding paragraph motivates the next problem.

**Problem 5.6.** *Assuming $P \neq NP$, is the following problem decidable? Given a finite algebra $\mathbf{A}$, is TER-EQ($(A; \mathcal{F})$) in P?*

We end with an example that indicates the wide separation between between polynomial equivalence problems and the term equivalence problems.

**Example 5.7.** *There exists an infinite sequence of finite semigroups $\mathbf{A_1} \leq \ldots \mathbf{A_k} \leq \mathbf{A_{k+1}} \ldots$ such that each semigroup in the sequence generates the same variety $\mathcal{V}$ (i.e. for all positive integers $k$, $\mathcal{V}(\mathbf{A_k})= \mathcal{V}(\mathbf{A_{k+1}})$), and for every positive integer $k$, POL-EQ($\mathbf{A_{2k}}$) is in P, while POL-EQ($\mathbf{A_{2k+1}}$) is co-NP-complete.*

Let $M$ be an $m \times n$ regular 0-1 matrix. Let $N$ be the regular $(m+1) \times (n+1)$ 0-1 matrix such that $N(m+1, k) = 1 = N(j, n+1)$ ($1 \leq k \leq n+1, 1 \leq j \leq m+1$) and $M$ is the upper left corner of $N$. We claim that POL-EQ($\mathbf{S_N}$) is in P.

We make an observation about concerning Z-sets for $\mathbf{S_N}$. Suppose $p$ is a polynomial over $\mathbf{S_N}$, with associated bipartite graph $B(p)= (X, Y, E)$. See the paragraph after Observation 4.1 for notation. Assume that $\tilde{f}$ is a partial function from $X \cup Y \to I \cup \Lambda$. We make an observation: There exists an extension $\tilde{e}$ of $\tilde{f}$, $\tilde{e} : X \cup Y \to I \cup \Lambda$ (so $\tilde{e}$ is a total function with domain $X \cup Y$) such that $e(p) \neq 0$ where $e$ is the evaluation determined by $\tilde{e}$, if and only if the function $\tilde{f}' : X \cup Y \to I \cup \Lambda$ which extends $\tilde{f}$ where $\tilde{f}'(x) = n + 1$, if $x$ is not in the domain of $\tilde{f}'$, and $\tilde{f}'(y) = m+1$, if $y$ is not in the domain of $\tilde{f}$, satisfies $f'p) \neq 0$, where $f'$ is the evaluation determined by $\tilde{f}'$.

It is easy to see that POL=0($\mathbf{S_N}$) is in P. Given a polynomial $p$ over $\mathbf{S_M}$, $p$ is identically 0 if and only if the evaluation $e$ which such that for each variable $v$ occuring in $p$, $e(v) = [n+1, m+1]$, has the property that $e(p) \neq 0$. In fact the evaluation $e$ is a a special case of the scenario described in the previous paragraph, one in which the partial function $\tilde{f}$ has domain the empty set. So to check if $p$ is identically 0 it suffices to check if a single evaluation $e$ satisfies $e(p) \neq 0$. This can be done in $length(p) - 1$ steps. It follows that POL=0($\mathbf{S_N}$) is in P.

Let $p, q$ over $\mathbf{S_N}$. Since POL=0($\mathbf{S_N}$) is in P, without loss of generality we may assume that neither $p$ nor $q$ is identically 0. If we assume that $p$ nor $q$ is identically 0, then we may assume also that the variables occuring in $p$ and in $q$ are the same (if not, it is easy to see that $Z(p) \neq Z(q)$).

Let $a, b$ be a pair of non-edge vertices of $B(p)$ with $a \in X$ and $b \in Y$. If there exists an evaluation $e$ of the variables of $p$ such that $M(\tilde{e}(a), \tilde{e}(b)) = 0$ and $e(p) \neq 0$, then we say that $a, b$ are *p-matchable*. Suppose that there exists a pair of non-edge vertices $a, b$ of $B(p)$ where $a \in X$, $b \in Y$, such that $a, b$ is $p$-matchable and $a, b$ **is** an edge of $B(q)$. Then we claim that $Z(q)$ is not a subset of $Z(p)$. To see this, let $e$ be an evaluation of the variables of $p, q$ such that $M(\tilde{e}(a), \tilde{e}(b)) = 0$ and $e(p) \neq 0$. Since $a, b$ is an edge of $B(q)$ and $M(\tilde{e}(a), \tilde{e}(b)) = 0$, it follows that $e(q) = 0$. This completes the proof of the claim. Consider the following property of the graphs $B(p)$ and $B(q)$.

(\*\*)If $a, b$ is an edge of $B(q)$ but $a, b$ is a pair of non-edge vertices of $B(p)$, $a \in X \subset B(p)$, $b \in Y \subset B(p)$, then $a, b$ is not $p$-matchable.

If $B(p)$ and $B(q)$ satisfy (\*\*), then we claim that $Z(q) \subset Z(p)$. If $e(q) = 0$, then there exist $a, b$ an edge of $Z(q)$ such that $M(\tilde{e}(a), \tilde{e}(b)) = 0$. Since $a, b$ is not $p$-matchable, it follows that $e(p) = 0$. As we saw in the paragraph preceding (\*\*), if neither $p$ nor $q$ is identically 0, and (\*\*) is not satisfied, then $Z(q)$ is not a subset of $Z(p)$.

We claim that determining the non-edge pairs $B(p)$ that are $p$-matchable can be computed in polynomial time (with respect to $length(p)$). Let $a \in X \subset B(p)$, $b \in Y \subset B(p)$ and let $a, b$ be a non-edge of $B(p)$. Let $f$ be a partial map of $X \cup Y$ with domain $a, b$ such that $M(\tilde{f}(b), \tilde{f}(a)) = 0$ (note that there are at most $mn$ such partial maps). We can extend $\tilde{f}$ to a total map $\tilde{e}$ with domain on $X \cup Y$ if and only if the evaluation $f'$ associated with $\tilde{f}' : X \cup Y \to I \cup \Lambda$, where $\tilde{f}'(x) = n + 1$, if $x$ is not in the domain of $\tilde{f}'$, and $\tilde{f}'(y) = m + 1$, if $y$ is not in the domain of $\tilde{f}$, satisfies $f'(p) \neq 0$. In particular, to check whether $a, b$ is $p$-matchable involves at most $mn\ length(p-1)$ computations. It now follows we can determine if (\*\*) is satisfied in polynomial time.

We can determine if $p$, $q$ are identically 0, determine the graphs $B(p)$, $B(q)$, the $p$-matchable non-edge pairs of $B(p)$, the $q$-matchable non-edge pairs of $B(q)$, and compare these to determine whether (\*\*) is satisfied, and whether (\*\*) is satisfied with the roles of $p$ and $q$ switched, all in polynomial time. By the preceding paragraphs, we can conclude that POL=Z-SET($\mathbf{S_N}$) is in P. From the first paragraph of the statement of Lemma 4.13, it follows that POL-EQ(($\mathbf{S_N}$) is in P.

Let $\mathbf{A_1} = \mathbf{S_{H_3}}$. Let $\mathbf{A_2} = \mathbf{S_N}$ where $N$ is the $H_3$ matrix, "bordered with 1's" (so $N$ is a $4 \times 4$ matrix). Let $\mathbf{A_3}$ be the combinatorial Rees matrix semigroup $\mathbf{S_{N_3}}$ where $N_3$ is the direct sum of $N$ and $H_3$. We leave it to the reader to verify that POL-EQ($\mathbf{S_{N_3}}$) is co-NP-complete ( in fact, GRAPH 3-COLORABILITY polynomially embeds into POL=0($\mathbf{S_{N_3}}$) by mimicking the proof of Lemma 3.1 but in the setting of $\mathbf{S_{N_3}}$). Let $\mathbf{A_4}$ be the combinatorial Rees matrix semigroup $\mathbf{S_{N_4}}$ where $N_4$ is the result of bordering $N_3$ with 1's; $\mathbf{A_5}$ is the combinatorial Rees matrix semigroup with a matrix that is the direct sum of $N_4$ with $H_3$, ad infinitum. We have for any positive integer $k$, POL-EQ($\mathbf{A_{2k}}$) is in P while POL-EQ($\mathbf{A_{2k+1}}$) is co-NP-complete.

A semigroup $\mathbf{S}$ is said to be *congruence-free* if the only congruences of $\mathbf{S}$ are the diagonal and universal relations. Universal algebraists call such semigroups (algebras) *simple*. It is not difficult to show that for every positive integer $k$, $\mathbf{A_k}$ above, is congruence-free (as is any combinatorial Rees matrix semigroup $\mathbf{S_M}$ as long as no two rows of $M$ are equal and no two columns of $M$ are equal).

**Problem 5.8.** *Find conditions on finite simple algebras which guarantee that whenever POL-EQ ($\mathbf{A}$) is in P and $\mathbf{B}$ is a subalgebra of $\mathbf{A}$, then POL-EQ ($\mathbf{B}$) is in P.*

As we remarked after Theorem 1.1 in Section 1.2, the complexity of polynomial equivalence for a finite ring depends entirely on nilpotence, a property which is inherited by subrings. Thus polynomial equivalence for finite rings is nicely behaved with respect to subrings. We assume some familiarity with basic universal algebra for parts of the next problem. It is not difficult to see that t e variety generated by any of the algebras of Example 5.7 is not congruence permutable. The first part of the next problem is a broadening of Problem 5.8.

**Problem 5.9.** *Find conditions on a variety $\mathcal{V}$ so that whenever $\mathbf{T}$ is a finite algebra in $\mathcal{V}$ and $\mathbf{S}$ is a subalgebra of $\mathbf{T}$, then POL-EQ($\mathbf{T}$) is in P implies that POL-EQ($\mathbf{S}$) is in P. Is it true that this property holds in congruence permutable varieties?*